\documentclass[10pt,twocolumn]{article}

\usepackage[utf8]{inputenc}
\usepackage[T1]{fontenc}
\usepackage{amsmath,amssymb,amsthm,bm}
\usepackage{graphicx}
\usepackage{hyperref}
\usepackage{xcolor}
\usepackage{listings}
\usepackage{geometry}
\usepackage{booktabs,array,multirow,makecell}
\usepackage{enumitem}
\usepackage[numbers]{natbib}
\usepackage{titlesec}
\usepackage{parskip}
\usepackage{caption}
\usepackage{float}
\usepackage{comment}
\usepackage{quantikz}
\usetikzlibrary{positioning,arrows.meta}
\usepackage{flafter} 
\usepackage[section]{placeins}
\titlespacing*{\section}{0pt}{8pt}{4pt}
\titlespacing*{\subsection}{0pt}{4pt}{3pt}
\titlespacing*{\subsubsection}{0pt}{3pt}{3pt}
\titlespacing*{\paragraph}{0pt}{2pt}{1pt}
\newcolumntype{M}{>{$}l<{$}}

\definecolor{codeblue}{rgb}{0.13,0.29,0.53}
\definecolor{codegray}{rgb}{0.5,0.5,0.5}
\definecolor{codegreen}{rgb}{0.0,0.45,0.0}
\definecolor{codeback}{rgb}{0.97,0.97,0.97}

\hypersetup{colorlinks=true, linkcolor=codeblue,
            citecolor=codegreen, urlcolor=codeblue}

\newcommand{\pkgname}{QuantumPoissonHomogenization}
\newcommand{\libname}{quantumpoissonhomogenization}
\newcommand{\repourl}{https://github.com/UW-ERSL/QuantumPoissonHomogenization}

\newcommand{\kron}{\otimes}
\newcommand{\norm}[1]{\ensuremath{\|#1\|}}
\newcommand{\abs}[1]{\ensuremath{|#1|}}
\newcommand{\mat}[1]{\mathbf{#1}}
\newcommand{\Kchi}{\mat{K}^{\chi}}

\newtheorem{proposition}{Proposition}

\theoremstyle{definition}
\newtheorem{example}{Example}

\title{\textbf{Mesh-Independent Preconditioning of Quantum Linear Solvers
for the Periodic Two-Phase Poisson Problem}}
\author{Krishnan Suresh\\
\small University of Wisconsin--Madison\\
\small \texttt{ksuresh@wisc.edu}}
\date{}

\begin{document}
\maketitle
\thispagestyle{empty}

\begin{abstract}
We consider the quantum solution of a periodic two-phase Poisson problem.
Quantum linear-system algorithms for its finite element
discretization are limited by the condition number of the stiffness matrix.
This condition number grows quadratically with the number of elements per
side. Classically, preconditioning by the homogeneous (single-phase) operator bounds the
condition number by the conductivity contrast, for every mesh size. On a
quantum computer, however, a preconditioner encoded as a separate factor
multiplies the subnormalizations, and the quadratic dependency returns.

In this paper, we absorb the closed-form inverse square root of the
homogeneous operator into the finite element factors. The preconditioned operator then becomes the product of a norm-preserving gradient operator $\mat{V}$, a diagonal matrix of scaled conductivities, and
$\mat{V}^T$. Each factor has norm at most one. The block encoding therefore has subnormalization one, and its
condition number is bounded by just the conductivity contrast, for every mesh size and in every
dimension. The operator $\mat{V}$ is applied as a quantum Fourier transform, a small
rotation at each wavenumber, and an inverse transform. Quantum singular value transformation
then block-encodes the inverse of the preconditioned operator.  

We test the method numerically on three periodic cells with up to $10^6$ unknowns.
The results confirm that the condition number does not depend on the mesh, and the
energy-norm error of the inverse is below the polynomial tolerance at every
mesh. On meshes of $4 \times 4$ and $8 \times 8$ elements, we also build the
entire inverse as one Qiskit circuit and simulate it by state vector, with the
rotation at each wavenumber built from uniformly controlled rotations. The
output state agrees with the finite element solution to within the polynomial
tolerance. 
\end{abstract}


\section{Introduction}
\label{sec:intro}

Quantum linear-system algorithms prepare a quantum state proportional to the
solution of a linear system $\mat{A}\bm x = \bm b$. Following the work of
Harrow, Hassidim, and Lloyd~\citep{harrow2009quantum}, later algorithms have
improved the dependence of the cost on the condition number and on the
precision \citep{childs2017quantum,gilyen2019quantum,martyn2021grand,
costa2022optimal,tosti2022review}. In their modern form, the matrix is
supplied as a \emph{block encoding}, a unitary circuit that contains a scaled
copy of $\mat{A}$ as one of its blocks, and the quantum singular value
transformation (QSVT) applies a polynomial approximation of $1/x$ to it. If the condition number of $\mat{A}$ is
bounded and its block encoding is efficient, the cost of these algorithms grows only polylogarithmically with the
number of unknowns \citep{harrow2009quantum,childs2017quantum}.
Linear systems from discretized partial differential equations are natural
candidates for these algorithms, since their size grows rapidly with mesh
refinement \citep{montanaro2016quantum,liu2024towards}.
For such systems, Aaronson~\citep{aaronson2015read} lists four obstacles that can
remove the polylogarithmic advantage end to end: loading the right-hand side, encoding the operator, the
condition number, and reading out the solution. This paper is concerned with
the third: the condition number. The first, loading the right-hand side, is
addressed for structured vectors in our earlier work~\citep{suresh2026pyencode},
and the fourth, extracting a quantity of interest from the output state, is the
subject of future work.

For the stiffness matrix of a Poisson problem on a mesh with
$N$ elements per side, the condition number grows as $\Theta(N^2)$. The number of queries to the block encoding grows
at least linearly with the condition number, so it grows at least in
proportion to $N^2$. The classical remedy is preconditioning: a matrix that is
cheap to apply and spectrally close to the inverse is multiplied onto the
system, and the condition number of the product no longer depends on $N$. On
a quantum computer, the preconditioner, if applied as a separate factor, must itself be supplied as a block
encoding, and the subnormalizations of the block encodings multiply. For several
natural preconditioners, the subnormalization of the preconditioner grows with
the mesh. The product of the subnormalizations then raises the condition number
seen by the circuit back to at least $\Theta(N^2)$, and the preconditioned
solve needs at least as many queries as the unpreconditioned one
\citep{golden2022quantum,lapworth2025preconditioned,pechan2026block}. Whether
a preconditioner helps on a quantum computer therefore depends on how it
enters the block encoding, not only on its classical properties.

In this paper, we ask whether the homogeneous (single-phase) preconditioner can be 
carried into a quantum block encoding without the
multiplication of subnormalizations, and what the resulting circuit costs when it
is built from gates. We consider periodic two-phase Poisson problem on a structured mesh of bilinear elements, and construct the block encoding of the preconditioned operator
and of its inverse.

In Section~\ref{sec:background}, we set up the two-phase problem, its
classical preconditioner, and the block-encoding and QSVT baseline, and show
why a separately encoded preconditioner fails. In
Section~\ref{sec:proposed}, we factor the
preconditioned operator as $\mat{V}^T\tilde{\mat{\Sigma}}\,\mat{V}$, with
every factor of norm one, and derive the Fourier structure of $\mat{V}$ in
one and two dimensions. We then construct its circuit, block-encode the
operator and its inverse, and give the gate-level circuit, its precision, and
its cost; the wavenumber-controlled stage is constructed in
Appendix~\ref{app:multiplexed}. The open-source implementation and its use
are described in Section~\ref{sec:implementation}. Numerical experiments on
three cells with up to $10^6$ unknowns, and a state-vector simulation of the
entire circuit on small meshes, are reported in
Section~\ref{sec:numerics}, and the method is assessed in
Section~\ref{sec:conclusions}. Table~\ref{tab:symbols} lists the principal
symbols used in this paper.

\section{Problem Formulation and Obstacles}
\label{sec:background}

In this section, we formulate the two-phase problem and its preconditioned
operator and examine why the classical construction does not directly carry
over to block encodings.

\begin{table*}[t]
\caption{Principal symbols. Bold upright letters denote matrices and
vectors, calligraphic letters denote unitaries and the Fourier transform,
and $(\cdot)^{+}$ denotes the pseudo-inverse. The left columns cover the
two-phase problem and its finite element matrices; the right columns cover
the Fourier structure, the block encodings, and the inverse. Symbols used in
one subsection only are defined there and not listed, and the symbols of the
fixed-point arithmetic are defined in Appendix~\ref{app:multiplexed}.}
\label{tab:symbols}
\centering
\footnotesize
\renewcommand{\arraystretch}{1.15}
\begin{tabular}{@{}Mp{0.35\linewidth}Mp{0.35\linewidth}@{}}
\toprule
\text{Symbol} & Meaning & \text{Symbol} & Meaning\\
\midrule
\sigma_{\min},\ \sigma_{\max} & conductivities of matrix and inclusion &
\mathcal{F} & two-dimensional discrete Fourier transform, $N^2 \times N^2$\\
\rho & contrast $\sigma_{\max}/\sigma_{\min}$ &
\bm k & wavenumber\\
v_f & inclusion volume fraction &
\hat{S}(\bm k) & eigenvalue of $\mat{K}$ at wavenumber $\bm k$\\
\chi & phase indicator, $1$ in the inclusion &
\mat{V} & partial isometry $\mat{W}^{1/2}\mat{G}(\mat{K}^{+})^{1/2}$, $8N^2 \times N^2$\\
N,\ h & elements per side; $h = 1/N$ &
\hat{\mat{V}}(\bm k) & unit vector of length $8$ that forms $\mat{V}$\\
N_{\mathrm{dof}},\ m & unknowns $N^2$; qubits $\log_2 N_{\mathrm{dof}}$ &
\mathcal{Q}(\bm k) & $16 \times 16$ unitary with first column $\hat{\mat{V}}(\bm k)$\\
d & spatial dimension &
\mat{A} & $\Kchi/\sigma_{\max}$\\
\bm f & heat source, orthogonal to the constant fields &
\alpha & subnormalization of an encoding of $\mat{A}$\\
\Kchi,\ \mat{K} & two-phase and unit-conductivity stiffness matrices &
\kappa & condition number seen by QSVT\\
\mat{G} & nodal values to gradients at quadrature points, $8N^2 \times N^2$ &
R_\chi & material reflection $I - 2\,\mathrm{diag}(\chi)$\\
\mat{W} & quadrature weights $w\,\mat{I}$, $w = h^2/4$ &
\mathcal{U}_V,\ \mathcal{U}_\Sigma,\ \mathcal{U}_M & block encodings of $\mat{V}$, $\tilde{\mat{\Sigma}}$, $\tilde{\mat{M}}$\\
\mat{\Sigma},\ \tilde{\mat{\Sigma}} & conductivity at quadrature points; $\mat{\Sigma}/\sigma_{\max}$ &
q_t & odd QSVT polynomial that approximates $1/x$\\
\bar{\sigma} & average of $\tilde{\mat{\Sigma}}$ over the cell, equal to \mbox{$v_f + (1-v_f)/\rho$} &
d_t,\ C & degree of $q_t$ (queries to $\mathcal{U}_M$); its subnormalization $\max_{[-1,1]}\abs{q_t}$\\
\mat{M},\ \tilde{\mat{M}} & preconditioned operator; its encoded form &
\epsilon_p & residual tolerance of $q_t$\\
\mat{\Pi},\ \mat{\Pi}_0 & projectors onto the complement of the constant fields, and onto them &
F,\ \delta & fractional bits; encoding error of $\mathcal{U}_V$\\
\bottomrule
\end{tabular}
\end{table*}

\subsection{Periodic Two-Phase Poisson Problem}
\label{sec:computational-homogenization}
\begin{figure}[H]
\centering
\begin{minipage}[b]{0.48\linewidth}
\centering
\includegraphics[width=0.42\linewidth]{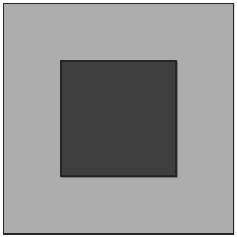}\\[2pt]
(a)
\end{minipage}
\begin{minipage}[b]{0.48\linewidth}
\centering
\includegraphics[width=0.42\linewidth]{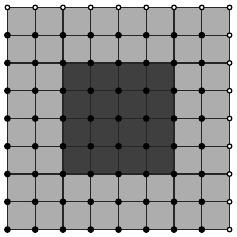}\\[2pt]
(b)
\end{minipage}
\caption{Worked example for this paper. (a) A square inclusion of volume
fraction $0.25$, centered in the unit cell. (b) The same cell discretized by a
structured $N \times N$ mesh of bilinear Q4 elements for $N = 4$; each element is assigned
to the inclusion or the matrix by the phase indicator $\chi$, and all elements
share one element matrix up to the scale $\sigma_{\min}$ or $\sigma_{\max}$.}
\label{fig:cell}
\end{figure}

Consider a periodic cell of two isotropic phases, an inclusion in a matrix, and the heterogeneous Poisson equation
$-\nabla\cdot(\sigma\nabla u) = f$ on the cell
(Figure~\ref{fig:cell}(a)). In engineering, this equation models
steady heat conduction, with $u$ the temperature, $\sigma$ the thermal
conductivity, and $f$ the heat source. The phases have conductivities
$\sigma_{\max}$ (inclusion) and $\sigma_{\min}$ (matrix), with contrast
$\rho = \sigma_{\max}/\sigma_{\min} \ge 1$; $v_f$ denotes the inclusion
volume fraction. Such cells are the unit cells of computational
homogenization \citep{hassani1998review,lucarini2021fft,
andreassen2014composite}, which computes their effective conductivity from
one cell problem per macroscopic gradient; this paper solves the cell problem
and does not compute the effective conductivity. The worked example of this paper is
a square inclusion of $0.25$ volume fraction, centered in the cell (Figure~\ref{fig:cell}(a))
that can be encoded at $O(1)$ cost~\citep{suresh2026quantum}.

A standard approach to solving the two-phase problem is to discretize the cell with a
structured $N \times N$ mesh of bilinear quadrilateral (Q4) elements, as in
Figure~\ref{fig:cell}(b), and perform a finite element analysis
(FEA). With one unknown per node, the periodic mesh has
$N_{\mathrm{dof}} = N^2$ degrees of freedom; with $N$ a power of two, these
are represented on $m = \log_2 N_{\mathrm{dof}}$ qubits. We denote the
periodic two-phase finite-element stiffness matrix by $\Kchi$.

\subsection{Classical Preconditioning}
\label{sec:reference-preconditioning}

The two-phase stiffness matrix $\Kchi$ is ill-conditioned under mesh
refinement. On a mesh of width $h = 1/N$ in $d$ dimensions, the largest
eigenvalue of $\Kchi$ on the complement of the constant fields is
$\Theta(h^{d-2})$ and the smallest is $\Theta(h^{d})$ \cite{ern2006evaluation}. The
condition number of $\Kchi$ therefore grows as $\Theta(N^2)$ in every
dimension. The conductivity contrast contributes a further factor of up to $\rho$. One remedy is to precondition $\Kchi$ by
the stiffness matrix $\mat{K}$ of the homogeneous cell with unit
conductivity
\citep{moulinec1998numerical,nielsen2009preconditioning,ladecky2023optimal}.
$\mat{K}$ is the natural choice because of two properties, discussed next.

\subsubsection{$\mat{K}$ is diagonalizable by the Fourier transform}

First, unlike $\Kchi$, the matrix $\mat{K}$ is circulant. It is therefore
diagonalized by the unitary two-dimensional discrete Fourier transform
$\mathcal{F}$, with entries
$\mathcal{F}_{\bm k \bm j} = N^{-1} e^{-2\pi\mathrm{i}\,\bm k\cdot\bm j/N}$
for the mesh nodes $\bm j$ and the wavenumbers $\bm k$ in
$\{0,\dots,N-1\}^2$:
\begin{equation}
\label{eq:symbol}
\mat{K} = \mathcal{F}^{\dagger}\,
\mathrm{diag}_{\bm k}\bigl(\hat{S}(\bm k)\bigr)\,
\mathcal{F},
\end{equation}
where $\hat{S}(\bm k)$, the eigenvalue of $\mat{K}$ at wavenumber $\bm k$, is
a trigonometric function of $\bm k$ known in closed form
\citep{fix2008analysis,leuschner2018fourier,ladecky2023optimal}.
Equation~\eqref{eq:symbol} is worked out on the smallest mesh in
Example~\ref{ex:symbol}.

\begin{example}
\label{ex:symbol}
For $N = 2$ in two dimensions, order the four nodes and the four wavenumbers as
$(0,0), (0,1), (1,0), (1,1)$. The assembled stiffness matrix and the
transform are
\begin{equation}
\label{eq:ex-KF}
\begin{gathered}
\mat{K} = \frac13
\begin{bmatrix}
8 & -2 & -2 & -4\\
-2 & 8 & -4 & -2\\
-2 & -4 & 8 & -2\\
-4 & -2 & -2 & 8
\end{bmatrix},
\\
\mathcal{F} = \frac12
\begin{bmatrix}
1 & 1 & 1 & 1\\
1 & -1 & 1 & -1\\
1 & 1 & -1 & -1\\
1 & -1 & -1 & 1
\end{bmatrix} = \mathcal{F}^{\dagger},
\end{gathered}
\end{equation}
since $e^{-2\pi\mathrm{i}\,\bm k\cdot\bm j/2} = \pm 1$. The columns of
$\mathcal{F}^{\dagger}$ are the constant field, the fields that alternate
in $y$ and in $x$, and the checkerboard. Multiplying each by $\mat{K}$
gives the eigenvalues $\hat{S} = 0$, $4$, $4$, and $8/3$, and
\begin{equation}
\label{eq:ex-S}
\mathcal{F}^{\dagger}\,\mathrm{diag}\bigl(0,\, 4,\, 4,\, \tfrac83\bigr)\,
\mathcal{F} = \mat{K} .
\end{equation}
The zero belongs to the constant field, which $\mat{K}$ maps to zero.
The pseudo-inverse, denoted $(\cdot)^{+}$, inverts the three nonzero entries
and keeps the zero:
\begin{equation}
\label{eq:ex-Kplus}
\mat{K}^{+} = \mathcal{F}^{\dagger}\,
\mathrm{diag}\bigl(0,\, \tfrac14,\, \tfrac14,\, \tfrac38\bigr)\,\mathcal{F}
= \frac{1}{32}
\begin{bmatrix}
7 & -3 & -3 & -1\\
-3 & 7 & -1 & -3\\
-3 & -1 & 7 & -3\\
-1 & -3 & -3 & 7
\end{bmatrix},
\end{equation}
and $(\mat{K}^{+})^{1/2}$ is formed in the same way from
$\mathrm{diag}\bigl(0,\, \tfrac12,\, \tfrac12,\, \sqrt{3/8}\bigr)$.
\end{example}

The same construction, including the zero eigenvalue, holds for every $N$.
The constant fields form the null
space of both $\Kchi$ and $\mat{K}$ and carry the zero eigenvalue of
$\mat{K}$; the smallest nonzero eigenvalue decreases as $h^2$
\citep{fix2008analysis}, so $\norm{\mat{K}^{+}}_2 = \Theta(N^2)$. By
Equation~\eqref{eq:symbol}, any function $f$ of $\mat{K}$ is
$\mathcal{F}^{\dagger}\,\mathrm{diag}_{\bm k}\bigl(f(\hat{S}(\bm k))\bigr)
\,\mathcal{F}$, with $f(0)$ replaced by $0$ \citep{suresh2026appliedqc}. In
particular, $\mat{K}^{+}$ and $(\mat{K}^{+})^{1/2}$, the matrix used in the
preconditioner below, are applied to a vector with two fast Fourier
transforms, at a cost of $O(N_{\mathrm{dof}} \log N_{\mathrm{dof}})$
operations.

On a quantum computer, the two transforms are quantum Fourier transforms, and
the cost lies in the diagonal factor. A matrix that is diagonalized by the
quantum Fourier transform and whose eigenvalues are known in closed form is
\emph{fast-invertible} in the sense of Tong et al.~\citep{tong2021fast}: a
block encoding of its inverse can be built at a cost independent of its
condition number \cite[Prop.~7]{tong2021fast}, and $\mat{K}$ is of this kind.
For general $N$, $\hat{S}(\bm k)$
is a trigonometric function of $\bm k$, and applying $\mat{K}^{+}$ or
$(\mat{K}^{+})^{1/2}$ requires evaluating its reciprocal or square root
coherently for all $N^2$ wavenumbers. This fixed-point arithmetic dominates
the gate count of the method (Section~\ref{sec:cost}).

\subsubsection{$\mat{K}$ and $\Kchi$ are spectrally equivalent}
The second property is that
$\mat{K}$ and $\Kchi$ are assembled from the same element matrix and differ
only in the conductivity assigned to each element, so $\mat{K}$ has the
same mesh-dependent spectrum as $\Kchi$, and preconditioning by $\mat{K}$
removes it. To establish this, observe that both admit the factored form
\citep{ladecky2023optimal}
\begin{equation}
\label{eq:factored}
\Kchi = \mat{G}^T \mat{W} \mat{\Sigma}\, \mat{G},
\qquad
\mat{K} = \mat{G}^T \mat{W} \mat{G},
\end{equation}
where $\mat{G}$ maps nodal values to gradients at the
quadrature points, $\mat{W} = w\,\mat{I}$ is the common quadrature weight
$w = h^2/4$, and $\mat{\Sigma}$ is diagonal with entry $\sigma_{\min}$ or
$\sigma_{\max}$ at each quadrature point according to the phase of its
element. With two gradient components at each of the four quadrature
points of each of the $N^2$ elements, $\mat{G}$ is $8N^2 \times N^2$, and
$\mat{W}$ and $\mat{\Sigma}$ are $8N^2 \times 8N^2$.

For a single element, with corners ordered $(0,0)$, $(1,0)$, $(1,1)$,
$(0,1)$ and local coordinates $(\xi, \eta) \in [-1,1]^2$, the gradients of
the four bilinear shape functions at a point $(\xi, \eta)$ form the
$2 \times 4$ matrix
\begin{equation}
\label{eq:B-matrix}
\mat{B}(\xi,\eta) = \frac{1}{2h}
\begin{bmatrix}
-(1-\eta) & 1-\eta & 1+\eta & -(1+\eta)\\
-(1-\xi) & -(1+\xi) & 1+\xi & 1-\xi
\end{bmatrix}.
\end{equation}
The element block of $\mat{G}$ stacks $\mat{B}$ at the four Gauss points
$(\xi_g, \eta_g) \in \{\pm 1/\sqrt3\}^2$, ordered as $(-,-)$, $(+,-)$,
$(+,+)$, $(-,+)$, into an $8 \times 4$ matrix $\mat{G}_e$, with the $x$- and
$y$-gradients at point $g$ in rows $2g-1$ and $2g$, and $\mat{G}$ places these
blocks in the rows of each element and the columns of its four nodes. With
$w = h^2/4$, the element stiffness matrix is
\begin{equation}
\label{eq:Ke}
\mat{G}_e^T\,w\,\mat{G}_e = \frac16
\begin{bmatrix}
4 & -1 & -2 & -1\\
-1 & 4 & -1 & -2\\
-2 & -1 & 4 & -1\\
-1 & -2 & -1 & 4
\end{bmatrix},
\end{equation}
independent of $h$, and Equation~\eqref{eq:factored} is the usual
assembly of these element matrices, each scaled by the conductivity of its
element.

Since $\mat{\Sigma}$ acts as a scalar at each quadrature point,  $\bm u^T \Kchi \bm u$ is a weighted sum over quadrature
points of $(\mat{G}\bm u)^T \mat{W} (\mat{G}\bm u)$ of the
unit-conductivity cell, and every nodal vector $\bm{u}$
satisfies
\begin{equation}
\label{eq:energy-bounds}
\sigma_{\min}\, \bm{u}^T \mat{K} \bm{u}
\;\le\; \bm{u}^T \Kchi \bm{u}
\;\le\; \sigma_{\max}\, \bm{u}^T \mat{K} \bm{u} .
\end{equation}
That is, the energy of the two-phase cell lies between those of a cell
made entirely of matrix and a cell made entirely of inclusion phase.

To carry these energy bounds over to a preconditioned operator, we apply
the preconditioner symmetrically, on both sides, so that the preconditioned
operator
\begin{equation}
\label{eq:M}
\mat{M} = \frac{1}{\sigma_{\max}}\,(\mat{K}^{+})^{1/2}\, \Kchi\,
(\mat{K}^{+})^{1/2}
\end{equation}
remains symmetric positive semidefinite. For $\bm v$ orthogonal to the
constant fields, set $\bm u = (\mat{K}^{+})^{1/2}\bm v$. Then
$\bm u^T\mat{K}\bm u = \bm v^T\bm v$ and $\bm v^T\mat{M}\bm v =
\bm u^T\Kchi\bm u / \sigma_{\max}$, and Equation~\eqref{eq:energy-bounds}
becomes
\begin{equation}
\label{eq:rayleigh}
\frac{1}{\rho}
\;\le\; \frac{\bm v^T \mat{M} \bm v}{\bm v^T \bm v}
\;\le\; 1 .
\end{equation}
The Rayleigh quotient of $\mat{M}$ is thus bounded by the contrast alone,
so the nonzero spectrum of $\mat{M}$ lies in $[1/\rho,\, 1]$ for every
mesh size $N$ and every phase arrangement $\chi$. Thus, classically, preconditioning
removes the mesh dependence.

The bound of Equation~\eqref{eq:rayleigh} is attained at both ends.
Reversing the substitution above, take $\bm u$ with $\bm v = \mat{K}^{1/2}\bm u$.
If the gradient of $\bm u$ vanishes in every inclusion element, those elements
contribute nothing to either $\Kchi\bm u$ or $\mat{K}\bm u$, the remaining
terms of the assembly differ only by the factor $\sigma_{\min}$, and
$\Kchi\bm u = \sigma_{\min}\mat{K}\bm u$; then $\bm v$ is an eigenvector of
$\mat{M}$ with eigenvalue $1/\rho$. A field whose gradient vanishes in every
matrix element gives the eigenvalue $1$ in the same way. Only the fields that
vary across the interface give eigenvalues strictly between the ends.

These results are established classically
\citep{nielsen2009preconditioning,gergelits2019laplacian,ladecky2023optimal}.
They underlie the FFT-based homogenization solvers, in which
preconditioned conjugate gradients converges in
$O(\sqrt{\rho}\,\log(1/\epsilon))$ iterations independent of $N$, each
costing $O(N_{\mathrm{dof}} \log N_{\mathrm{dof}})$ operations
\citep{moulinec1998numerical,schneider2017fft,leuschner2018fourier,
ladecky2023optimal}.

\subsection{Quantum Block Encoding}
\label{sec:block-encoding}

Before we ask if the classical preconditioner carries over to quantum, we
must address a more basic requirement: a quantum computer applies only
unitary operations, and the stiffness matrix $\Kchi$ is not unitary. The
standard way to act with a non-unitary matrix is to embed it, suitably
scaled, as the leading block of a larger unitary. A unitary $\mathcal{U}$
acting on $a$ ancilla qubits and the $m$ system qubits is a \emph{block
encoding} of a matrix $\mat{A}$ with \emph{subnormalization} $\alpha$
\citep{gilyen2019quantum} if
\begin{equation}
\label{eq:be}
\mat{A} = \alpha\,
\bigl(\langle 0|^{\otimes a}\otimes I\bigr)\,\mathcal{U}\,
\bigl(|0\rangle^{\otimes a}\otimes I\bigr),
\end{equation}
that is,
\begin{equation}
\label{eq:be-block}
\mathcal{U} =
\begin{bmatrix}
\mat{A}/\alpha & \ast\\
\ast & \ast
\end{bmatrix},
\end{equation}
where the blocks marked $\ast$ are whatever is needed to make
$\mathcal{U}$ unitary. 

Operationally, preparing the ancillas in $|0\rangle^{\otimes a}$ and the
system in a state $|\bm\phi\rangle$, applying $\mathcal{U}$, and measuring
the ancillas in $|0\rangle^{\otimes a}$ leaves the system in the state
proportional to $\mat{A}|\bm\phi\rangle$, with success probability
$\norm{\mat{A}\bm\phi}^2/\alpha^2$ \citep{suresh2026appliedqc}. Two consequences follow. First, since
every block of a unitary has spectral norm at most $1$, the
subnormalization satisfies $\alpha \ge \norm{\mat{A}}_2$ for every block
encoding, however it is built. Second,  a larger
$\alpha$ lowers the success probability, and, as shown in
Section~\ref{sec:qsvt}, it multiplies the condition number. A useful
block encoding therefore needs a circuit whose gate count grows at most
polylogarithmically with $N_{\mathrm{dof}}$ and a subnormalization close
to $\norm{\mat{A}}_2$.

For the periodic two-phase operator, an encoding that meets both
requirements follows the construction of~\citep{suresh2026quantum}:
$\Kchi/\sigma_{\max}$ is written as a linear combination of unitaries
built from cyclic shifts of the position register and the \emph{material
reflection} $R_\chi = I - 2\,\mathrm{diag}(\chi)$, which multiplies each
inclusion element by $-1$ and each matrix element by $+1$. The number of
terms is fixed by the Q4 stencil, and the subnormalization $\alpha$ of
this encoding of $\Kchi/\sigma_{\max}$, the meaning of $\alpha$ hereafter,
is independent of $N$ and bounded uniformly in $\rho$; the cell geometry
enters only through $R_\chi$, which for a square inclusion with dyadic
side length costs $O(1)$ gates~\citep{suresh2026quantum}. The details of
this encoding are not needed here; the cost of inverting it by QSVT is given
in Section~\ref{sec:qsvt}.

\subsection{Quantum Singular Value Transformation}
\label{sec:qsvt}

The quantum singular value transformation (QSVT)
\citep{gilyen2019quantum,martyn2021grand} applies a polynomial to a
block-encoded matrix. Given a block encoding $\mathcal{U}$ of a symmetric
matrix $\mat{A}$ with subnormalization $\alpha$, and a polynomial $p$ with
$|p(x)| \le 1$ on $[-1,1]$, QSVT uses a number of queries to $\mathcal{U}$ equal to the degree of $p$
to block-encode
\begin{equation}
\label{eq:qsvt-poly}
p(\mat{A}/\alpha) = \sum_j p(\lambda_j)\, \bm u_j \bm u_j^T ,
\end{equation}
where $\lambda_j$ and $\bm u_j$ are the eigenvalues and eigenvectors of
$\mat{A}/\alpha$. For $p(x) \approx 1/x$, the polynomial must follow $1/x$
down to the smallest eigenvalue of $\mat{A}/\alpha$, so the degree is set by
the condition number seen by the circuit,
\begin{equation}
\label{eq:kappa-circuit}
\kappa = \alpha\,\norm{\mat{A}^{+}}_2 ,
\end{equation}
as $d_t = O(\kappa \log(\kappa/\epsilon))$ for accuracy $\epsilon$
\citep{childs2017quantum, gilyen2019quantum}. We use the odd Chebyshev-iteration polynomial $q_t$ of degree $2t - 1$ of Gribling et al.~\cite[Lem.~13]{gribling2021improving}, which minimizes the
relative residual $\max_{1/\kappa \le |x| \le 1} |x\,q_t(x) - 1|$ among odd
polynomials of that degree.

For $\mat{A} = \Kchi/\sigma_{\max}$, $\alpha \ge \norm{\Kchi}_2/\sigma_{\max}$
(Section~\ref{sec:block-encoding}), so $\kappa$ is at least the condition
number of $\Kchi$, which grows as $\Theta(N^2)$
(Section~\ref{sec:reference-preconditioning}). For the encoding
of~\citep{suresh2026quantum}, $\alpha$ exceeds $\norm{\Kchi}_2/\sigma_{\max}$
by a factor bounded in $N$, so $\kappa$, and with it $d_t$, grows as
$\Theta(N^2)$.

\subsection{Quantum Preconditioning}
\label{sec:quantum-preconditioning}

Multiplying block encodings multiplies their subnormalizations \cite{gilyen2019quantum}, so a
preconditioner encoded as a separate factor cannot lower the condition
number seen by QSVT. Under block
encoding, $(\mat{K}^{+})^{1/2}$ must itself be encoded as the leading block
of a unitary, with a subnormalization
$\alpha_P \ge \norm{(\mat{K}^{+})^{1/2}}_2 = \norm{\mat{K}^{+}}_2^{1/2}$
(Section~\ref{sec:block-encoding}), which grows as $N$. In this subsection, we
show that this scaling cancels the benefit of the preconditioner and relate the
result to prior work.

\paragraph{The product barrier.}
Let $\mathcal{U}$ and $\mathcal{U}_P$ block-encode
$\mat{A} = \Kchi/\sigma_{\max}$ and $(\mat{K}^{+})^{1/2}$ with
subnormalizations $\alpha$ and $\alpha_P$. The product
$\mathcal{U}_P\,\mathcal{U}\,\mathcal{U}_P$ block-encodes
$(\mat{K}^{+})^{1/2}\mat{A}(\mat{K}^{+})^{1/2} = \mat{M}$
(Equation~\eqref{eq:M}) with subnormalization $\alpha_P^2\,\alpha$
\citep{gilyen2019quantum}, so the condition number seen by QSVT
(Equation~\eqref{eq:kappa-circuit}) is
$\kappa_P = \alpha_P^2\,\alpha\,\norm{\mat{M}^{+}}_2$. On the complement of
the constant fields, $\mat{A}^{+} = (\mat{K}^{+})^{1/2}\,\mat{M}^{+}
(\mat{K}^{+})^{1/2}$, so $\norm{\mat{A}^{+}}_2 \le \alpha_P^2\,
\norm{\mat{M}^{+}}_2$. Hence
\begin{equation}
\label{eq:product-barrier}
\kappa_P \;=\; \alpha_P^2\,\alpha\,\norm{\mat{M}^{+}}_2
\;\ge\; \alpha\,\norm{\mat{A}^{+}}_2 \;=\; \kappa :
\end{equation}
the preconditioned solve is never cheaper than the unpreconditioned one of
Section~\ref{sec:qsvt}, and $\kappa_P$ grows at least as $\Theta(N^2)$.
Pechan et al.~\citep{pechan2026block} prove the same inequality for every
preconditioner encoded as separate factors.

\paragraph{Prior work.}
The barrier applies to the following prior approaches.
\begin{enumerate}[leftmargin=1.4em,itemsep=2pt,topsep=2pt]
\item \textbf{Preconditioned block encodings:} Lapworth et
al.~\citep{lapworth2025preconditioned} report it numerically for
preconditioned quantum linear solvers.
\item \textbf{Inverse-Laplacian preconditioning:} Golden et
al.~\citep{golden2022quantum} precondition hydrological flow, a scalar
elliptic problem of the type treated here, by the inverse Laplacian, and
find that its norm grows linearly in the number of unknowns, so that its
block encoding carries the mesh dependence it removes classically.
\item \textbf{FFT-based fixed-point iteration:} Liu et
al.~\citep{liu2024towards} implement the iteration of Moulinec and
Suquet~\citep{moulinec1998numerical} for the cell problem with quantum
Fourier transforms, encoding each of its $S$ steps as a separate
non-unitary map. The $S$ steps form a polynomial of degree $S$ in the
preconditioned operator, composing them multiplies their
subnormalizations, and the success probability of the composed circuit is
not analyzed.
\item \textbf{Wavelet preconditioning:} Bagherimehrab et
al.~\citep{bagherimehrab2023fast} and Gharat et al.~\citep{gharat2026slac}
transform the operator to a wavelet basis and apply a diagonal
preconditioner, which bounds the condition number of the preconditioned
matrix independently of $N$. The cost of QSVT depends instead on the
condition number seen by the circuit, which includes the subnormalization of
the preconditioner and, for factors encoded separately, is bounded below by
Equation~\eqref{eq:product-barrier} \citep{pechan2026block}.
\item \textbf{Fast inversion:} for an operator $\mat{K}_0 + \mat{B}$ with
$\mat{K}_0$ fast-invertible, Tong et al.~\citep{tong2021fast} apply QSVT to
the preconditioned operator $I + \mat{K}_0^{-1}\mat{B}$, with a number of
queries that grows with the product of the subnormalizations of
$\mat{K}_0^{-1}$ and $\mat{B}$ \cite[Cor.~14]{tong2021fast}. This product does
not grow with the mesh when $\mat{B}$ is a lower-order term, such as a bounded
potential added to the Laplacian. For the two-phase operator, the natural
splitting is $\mat{K}_0 = \sigma_{\min}\mat{K}$, a cell made entirely of the
matrix phase, and $\mat{B} = \Kchi - \sigma_{\min}\mat{K}$, the additional
stiffness of the inclusion. The norm of $\mat{B}$ is proportional to the
largest eigenvalue of $\mat{K}$, and that of $\mat{K}_0^{+}$ to the inverse of
its smallest nonzero eigenvalue, so their product, and with it the query
count, grows as $\Theta(N^2)$.
\end{enumerate}

Deiml and Peterseim~\cite{deiml2026quantum} avoid the barrier by building the
preconditioner into the encoding, with a multilevel construction that does
not require periodicity, at a subnormalization that grows with the fine-scale
resolution and with the coefficient contrast. A preconditioner therefore lowers the cost of QSVT only if it is built into the block encoding of the preconditioned operator, not multiplied onto
an encoding of $\Kchi$; such an encoding is built in
Section~\ref{sec:proposed}.

\section{Fast-Inversion Block Encoding}
\label{sec:proposed}

Several concepts and equations were introduced in
Section~\ref{sec:background}, and four of them are used in this section.
\begin{itemize}[leftmargin=1.4em,itemsep=2pt,topsep=2pt]
\item \textbf{Factored form:} $\Kchi = \mat{G}^T\mat{W}\mat{\Sigma}\,\mat{G}$
and $\mat{K} = \mat{G}^T\mat{W}\mat{G}$ (Equation~\eqref{eq:factored}).
\item \textbf{Fourier diagonalization:} $\mat{K} = \mathcal{F}^{\dagger}\,
\mathrm{diag}_{\bm k}\bigl(\hat{S}(\bm k)\bigr)\,\mathcal{F}$
(Equation~\eqref{eq:symbol}), so functions of $\mat{K}$, such as
$(\mat{K}^{+})^{1/2}$, are diagonal in the Fourier basis.
\item \textbf{Preconditioned operator:} $\mat{M}$ of Equation~\eqref{eq:M},
whose nonzero eigenvalues lie in $[1/\rho, 1]$ for every $N$.
\item \textbf{Cost of QSVT:} the number of queries grows with the condition
number seen by the circuit, $\kappa = \alpha\,\norm{\mat{A}^{+}}_2$
(Equation~\eqref{eq:kappa-circuit}), which includes the subnormalization;
encoding $(\mat{K}^{+})^{1/2}$ as a separate factor makes it grow as $N^2$
(Equation~\eqref{eq:product-barrier}).
\end{itemize}

The aim of this section is a block encoding of $\mat{M}$ with
subnormalization $1$, so that the condition number seen by QSVT is at most
$\rho$ for every $N$, and a block encoding of its inverse. The construction
follows the FFT-based solvers of homogenization step by step, and
Table~\ref{tab:analogy} lists the correspondence. In those solvers, a
transform to wavenumbers, one small operation at each wavenumber, and a
transform back apply the preconditioner in $O(N_{\mathrm{dof}}\log
N_{\mathrm{dof}})$ operations. The circuit has the same three stages, with
quantum Fourier transforms in place of FFTs and a rotation of a small register
in place of the multiplication at each wavenumber. 

\begin{table}[H]
\caption{The FFT-based classical solver and the block-encoded circuit of this
section, step by step. Read each row across: the right column is the quantum
counterpart of the left. The costs are counts per application or per solve, not
run times, and the two columns are not comparable in absolute terms.}
\label{tab:analogy}
\centering
\footnotesize
\renewcommand{\arraystretch}{1.15}
\setlength{\tabcolsep}{3pt}
\begin{tabular}{@{}>{\raggedright\arraybackslash}p{0.21\linewidth}>{\raggedright\arraybackslash}p{0.35\linewidth}>{\raggedright\arraybackslash}p{0.38\linewidth}@{}}
\toprule
Step & FFT-based solver & Circuit (this paper)\\
\midrule
Preconditioner & $\mat{K}^{+}$ applied to a vector by two FFTs & $(\mat{K}^{+})^{1/2}$ absorbed into the gradient operator $\mat{V}$; never encoded on its own\\
To wavenumbers & FFT, $O(N_{\mathrm{dof}}\log N_{\mathrm{dof}})$ operations & quantum Fourier transform, $O(\log^2 N_{\mathrm{dof}})$ gates\\
At each wavenumber & multiply by $\hat{S}(\bm k)^{-1}$ & rotate the local register to the unit vector $\hat{\mat V}(\bm k)$\\
Back to the mesh & inverse FFT & inverse quantum Fourier transform\\
Conductivity & multiply by $\mat{\Sigma}$ at each quadrature point & $\mathcal{U}_\Sigma$, one controlled material reflection\\
Condition number & at most $\rho$, every $N$ & at most $\rho$, every $N$\\
Solve & conjugate gradients, $O(\sqrt\rho\log(1/\epsilon))$ iterations & QSVT, $O(\rho\log(1/\epsilon_p))$ queries\\
\bottomrule
\end{tabular}
\end{table}

In Section~\ref{sec:factorization}, we write $\mat{M}$ as a product
of three factors of norm $1$, $\mat{V}^T\tilde{\mat{\Sigma}}\,\mat{V}$. We
derive the Fourier structure of $\mat{V}$ in one and two dimensions in
Sections~\ref{sec:V-1D} and~\ref{sec:fourier-V}, and build its circuit in
Section~\ref{sec:circuit-V}. We encode $\tilde{\mat{\Sigma}}$ and $\mat{M}$ in
Section~\ref{sec:be-M} and the inverse in Section~\ref{sec:inverse}, and in
Section~\ref{sec:cost} we give the gate-level circuit for the
wavenumber-controlled stage, choose its precision, and state its cost.

\subsection{Unit-Norm Factorization}
\label{sec:factorization}

Our goal is a block encoding of the preconditioned operator $\mat{M}$ of
Equation~\eqref{eq:M} with subnormalization $1$. Unlike the fast inversion of
Tong et al.\ (Section~\ref{sec:quantum-preconditioning}), we do not
precondition from the left: we apply $(\mat{K}^{+})^{1/2}$ on both sides and
absorb it into the factors of $\Kchi$, so that the preconditioner is never
encoded as a separate factor. The idea is to apply the preconditioner to the
gradients rather than to the stiffness matrix.

We start from the factored form $\Kchi = \mat{G}^T\mat{W}\mat{\Sigma}\,\mat{G}$
of Equation~\eqref{eq:factored}. Substituting it into Equation~\eqref{eq:M}
gives
\begin{equation*}
\mat{M} = \frac{1}{\sigma_{\max}}\,(\mat{K}^{+})^{1/2}\,
\mat{G}^T\mat{W}\mat{\Sigma}\,\mat{G}\,(\mat{K}^{+})^{1/2} .
\end{equation*}
Two steps put this product in a symmetric form. First, we write
$\tilde{\mat{\Sigma}} = \mat{\Sigma}/\sigma_{\max}$ for the scaled
conductivity, a diagonal matrix whose entries are $1$ at the quadrature points
of the inclusion and $1/\rho$ at those of the matrix. Second, we split
$\mat{W} = \mat{W}^{1/2}\mat{W}^{1/2}$ and move one factor $\mat{W}^{1/2}$ to
the right of $\tilde{\mat{\Sigma}}$, which is allowed because
$\mat{W}^{1/2} = \sqrt{w}\,\mat{I}$ commutes with every matrix. The result is
\begin{equation}
\label{eq:M-regrouped}
\mat{M} =
\bigl[\mat{W}^{1/2}\mat{G}(\mat{K}^{+})^{1/2}\bigr]^T\,
\tilde{\mat{\Sigma}}\,
\bigl[\mat{W}^{1/2}\mat{G}(\mat{K}^{+})^{1/2}\bigr] .
\end{equation}
We denote the bracketed matrix by
\begin{equation}
\label{eq:V}
\mat{V} = \mat{W}^{1/2}\,\mat{G}\,(\mat{K}^{+})^{1/2} ,
\end{equation}
so that the preconditioned operator becomes
\begin{equation}
\label{eq:M-VSV}
\mat{M} = \mat{V}^T\,\tilde{\mat{\Sigma}}\,\mat{V} .
\end{equation}

Reading Equation~\eqref{eq:V} from right to left, $(\mat{K}^{+})^{1/2}$ maps
the $N^2$ nodal values to $N^2$ nodal values, $\mat{G}$ maps nodal values to
the $8$ gradient values of each element, and $\mat{W}^{1/2}$ multiplies each
gradient value by $\sqrt{w}$. The matrix $\mat{V}$ therefore has the size of
$\mat{G}$, $8N^2 \times N^2$. Its key property is that it preserves length.

\begin{proposition}[Unit-norm factorization]
\label{prop:factorization}
For every mesh size $N$ and every phase arrangement
$\chi$,
\begin{equation}
\label{eq:factorization}
\mat{V}^T\mat{V} = \mat{\Pi},
\qquad
\norm{\mat{V}}_2 = 1 ,
\end{equation}
where $\mat{\Pi}$ is the orthogonal projector onto the complement of the
constant fields.
\end{proposition}

\begin{proof}
From Equation~\eqref{eq:factored}, $\mat{V}^T\mat{V} =
(\mat{K}^{+})^{1/2}\mat{G}^T\mat{W}\mat{G}(\mat{K}^{+})^{1/2} =
(\mat{K}^{+})^{1/2}\mat{K}(\mat{K}^{+})^{1/2}$, which is the projector
onto the range of $\mat{K}$. Since $\mat{K}$ is symmetric, this range is
the orthogonal complement of its kernel, the constant fields. Finally, $\mat{V}^T\mat{V} = \mat{\Pi}$
means that the singular values of $\mat{V}$ are the square roots of the
eigenvalues of $\mat{\Pi}$, which are $0$ and $1$, so
$\norm{\mat{V}}_2 = 1$.
\end{proof}

In words, $\mat{V}$ preserves the Euclidean norm of every nodal vector
orthogonal to the constant fields and maps the constant fields to zero; a
matrix with this property is a \emph{partial isometry}
\citep{ben2003generalized}. Only Equation~\eqref{eq:factored} is used in the
proof, so the proposition holds in every dimension $d$.

By Equation~\eqref{eq:M-VSV}, $\mat{M}$ is thus a product of three factors
of norm $1$. The two copies of $\mat{V}$ depend on the mesh but not on the
microstructure. The diagonal factor $\tilde{\mat{\Sigma}}$ is the only factor
that depends on $\chi$, and its entries are $1$ and $1/\rho$; its subnormalization is therefore $1$. The
subnormalization of a product of block encodings is the product of their
subnormalizations, so block encodings of $\mat{V}$ and $\tilde{\mat{\Sigma}}$
with subnormalization $1$ give a block encoding of $\mat{M}$ with
subnormalization $1$. The encoding of $\mat{V}$ is constructed in
Sections~\ref{sec:V-1D} to~\ref{sec:circuit-V}, those of
$\tilde{\mat{\Sigma}}$ and $\mat{M}$ in Section~\ref{sec:be-M}, and that of
the inverse in Section~\ref{sec:inverse}. The gate-level circuit, its
precision, and its cost are then given in Section~\ref{sec:cost}.
\subsection{The Isometry in One Dimension}
\label{sec:V-1D}

Since $(\mat{K}^{+})^{1/2}$ is a dense matrix, $\mat{V}$ is dense, and a
circuit that loads its entries one by one would need a number of gates that
grows with the number of entries. What makes $\mat{V}$ cheap is its periodicity: $\mat{V}$ is circulant, so the Fourier transform
diagonalizes it, and the entire matrix is carried by one number per
wavenumber,
\begin{equation}
\label{eq:V-structure-1D}
\mat{V} = \mathcal{F}^{\dagger}\,
\mathrm{diag}_k\bigl(\hat{V}(k)\bigr)\,\mathcal{F} .
\end{equation}
This subsection derives $\hat{V}(k)$ in one dimension, where it is a phase.
Section~\ref{sec:fourier-V} repeats the derivation on the Q4 mesh, where the
diagonal becomes block diagonal.

Consider the periodic one-dimensional cell with $N$ linear elements of length
$h = 1/N$ and nodes $j = 0, \dots, N-1$ taken modulo $N$, with one quadrature
point per element of weight $w = h$. The gradient matrix is then
$N \times N$, with $(\mat{G}\bm u)_j = (u_{j+1} - u_j)/h$, so every row has
the same two entries shifted by one column from the row above. This is where
the periodicity enters: $\mat{G}$ is circulant, and so is every factor of
$\mat{V}$. Applying $\mat{G}$ to the plane wave $u_j = e^{\mathrm{i}\theta j}$,
with $\theta = 2\pi k/N$, multiplies it by
\begin{equation}
\label{eq:G-hat-1D}
\hat{G}(k) = \frac{e^{\mathrm{i}\theta} - 1}{h}
= \frac{2\mathrm{i}}{h}\, e^{\mathrm{i}\phi}\sin\phi ,
\qquad \phi = \frac{\theta}{2} = \frac{\pi k}{N} .
\end{equation}
The same plane wave is an eigenvector of $\mat{K}$, with eigenvalue
$\hat{S}(k) = w\,\abs{\hat{G}(k)}^2 = (4/h)\sin^2\phi$, and of
$\mat{W}^{1/2} = \sqrt{h}\,\mat{I}$. The three factors of $\mat{V}$ are
therefore diagonal in the same basis, and multiply wavenumber by wavenumber:
\begin{equation}
\label{eq:V-hat-1D}
\hat{V}(k) = \sqrt{h}\;\hat{G}(k)\;\hat{S}(k)^{-1/2}
= \mathrm{i}\,e^{\mathrm{i}\pi k/N},
\quad k = 1, \dots, N-1 ,
\end{equation}
since $\sin\phi > 0$ for these $k$, and $\hat{V}(0) = 0$ at the constant
field. Two things cancel: the element length $h$, so the result does not
depend on the mesh, and the magnitude, so every nonzero $\hat{V}(k)$ has
modulus $1$. The second is Proposition~\ref{prop:factorization} in this
setting, and it is what makes the subnormalization $1$: the inverse square
root $\hat S(k)^{-1/2}$, which alone grows as $N$, is absorbed exactly by the
growth of $\hat{G}(k)$.

The remaining structure is in the phase itself, which is linear in $k$ and
therefore factors over its bits: one single-qubit phase gate on each qubit of
the wavenumber register, and one more for the constant factor $\mathrm{i}$.
The one-dimensional $\mathcal{U}_V$ is thus Equation~\eqref{eq:V-structure-1D}
read as a circuit, a quantum Fourier transform, $n + 1$ phase gates with
$n = \log_2 N$, and the inverse transform, together with one gate, controlled
on $k = 0$, that moves the input to a state outside the encoded block
(Section~\ref{sec:circuit-V}). It needs no arithmetic.

In two dimensions, Equation~\eqref{eq:V-structure-1D} holds but the last
paragraph does not. An element has four nodes and eight gradient values, so
$\hat{V}(k)$ becomes a unit vector in eight dimensions and the diagonal
becomes block diagonal; and its entries depend on $\bm k$ through sines and
cosines of $\phi_1$ and $\phi_2$ that do not factor over the bits of
$\bm k$. We derive $\hat{\mat V}(\bm k)$ in Section~\ref{sec:fourier-V} and
prepare it with fixed-point arithmetic in Section~\ref{sec:cost}; this
arithmetic is the main cost of the entire method.

\subsection{The Isometry in Two Dimensions}
\label{sec:fourier-V}

We now repeat the steps of Section~\ref{sec:V-1D} for the Q4 mesh. Identify
each element with its lower-left node, so that element $\bm j$ has the nodes
$\bm j + \bm c$ for the four corner offsets $\bm c \in \{0,1\}^2$. The $8$
gradient values of element $\bm j$, two components at each of four Gauss
points, are linear combinations of the four nodal values $u_{\bm j + \bm c}$,
with coefficients from the matrix $\mat{B}$ of Equation~\eqref{eq:B-matrix}
that are the same for every element. As in one dimension, $\mat{G}$ is
therefore circulant, now in blocks: each element contributes a block of eight
rows, and the blocks repeat with a shift. A block-circulant matrix is
block-diagonalized by the discrete Fourier transform
\citep{davis1979circulant},
\begin{equation}
\label{eq:G-fourier}
\mat{G} = (\mathcal{F}^{\dagger}\kron I_8)\;
\mathrm{blockdiag}_{\bm k}\bigl(\hat{\mat{G}}(\bm k)\bigr)\;\mathcal{F}.
\end{equation}
where $I_8$ is the $8 \times 8$ identity matrix.

The block $\hat{\mat{G}}(\bm k)$ is found as in one dimension, by applying
$\mat{G}$ to the plane wave $u_{\bm j} = e^{\mathrm{i}\bm\theta\cdot\bm j}$:
each nodal value $u_{\bm j + \bm c}$ equals $e^{\mathrm{i}\bm\theta\cdot\bm c}$
times $u_{\bm j}$, and the eight gradient values of element $\bm j$ equal
$\hat{\mat{G}}(\bm k)$ times $u_{\bm j}$. With the half angles
$\phi_j = \pi k_j/N$ and the Gauss points
$(\xi_g, \eta_g) \in \{\pm 1/\sqrt3\}^2$, this gives the vector of length $8$
\begin{equation}
\label{eq:G-hat}
\hat{\mat{G}}(\bm k) = \frac{2\mathrm{i}}{h}\,e^{\mathrm{i}(\phi_1+\phi_2)}
\begin{bmatrix}
\sin\phi_1\,(\cos\phi_2 + \mathrm{i}\,\eta_g\sin\phi_2)\\
\sin\phi_2\,(\cos\phi_1 + \mathrm{i}\,\xi_g\sin\phi_1)
\end{bmatrix}_{g = 1,\dots,4},
\end{equation}
where the subscript $g$ denotes the four $2 \times 1$ blocks, one per Gauss
point, stacked into a vector of length $8$. The first entry of each block
is the $x$-gradient and the second the $y$-gradient. Both vanish at
$\bm k = \bm 0$, since a constant field has no gradient. The one-dimensional
result of Equation~\eqref{eq:G-hat-1D} has the same factor
$(2\mathrm{i}/h)\,e^{\mathrm{i}\phi}\sin\phi$.

The other two factors of $\mat{V}$ are already diagonal in the Fourier basis.
The weight $\mat{W}^{1/2} = \sqrt{w}\,\mat{I}$, with $w = h^2/4$, has the block
\begin{equation}
\label{eq:We}
\mat{W}_e^{1/2} = \sqrt{w}\,I_8 ,
\end{equation}
and $(\mat{K}^{+})^{1/2} = \mathcal{F}^{\dagger}\,
\mathrm{diag}_{\bm k}\bigl(\hat{S}(\bm k)^{-1/2}\bigr)\,\mathcal{F}$ by
Section~\ref{sec:reference-preconditioning}, with the entry for
$\bm k = \bm 0$ set to $0$. In the product
$\mat{V} = \mat{W}^{1/2}\mat{G}(\mat{K}^{+})^{1/2}$, the transforms in the
middle cancel, since $\mathcal{F}\mathcal{F}^{\dagger} = I$, and the
block-diagonal matrices multiply block by block. For every mesh size $N$,
this gives
\begin{equation}
\label{eq:V-fourier}
\begin{gathered}
\mat{V} = (\mathcal{F}^{\dagger}\kron I_{8})\;
  \mathrm{blockdiag}_{\bm k}\bigl(\hat{\mat{V}}(\bm k)\bigr)\;
  \mathcal{F},
\\
\hat{\mat{V}}(\bm k) = \mat{W}_e^{1/2}\,\hat{\mat{G}}(\bm k)\,
  \hat{S}(\bm k)^{-1/2},
\end{gathered}
\end{equation}
with $\hat{\mat{V}}(\bm 0) = \mat 0$, the two-dimensional counterpart of
Equation~\eqref{eq:V-hat-1D}. In one dimension, $\hat V(k)$ is a phase; here it
is a vector of length $8$, and it has unit length for every
$\bm k \ne \bm 0$. To see this, note that Equation~\eqref{eq:V-fourier} gives
$\mat{V}^T\mat{V} = \mathcal{F}^{\dagger}\,\mathrm{diag}_{\bm k}\bigl(
\norm{\hat{\mat{V}}(\bm k)}_2^2\bigr)\,\mathcal{F}$, since $\mat{V}$ is real.
By Proposition~\ref{prop:factorization}, this equals $\mat{\Pi}$, the same
product with $1$ on the diagonal except at the constant field $\bm k = \bm 0$.
Hence $\norm{\hat{\mat{V}}(\bm k)}_2 = 1$ for every $\bm k \ne \bm 0$: the
isometry is a unit vector at each wavenumber.

In explicit form, $\sqrt{w} = h/2$ cancels the factor $2/h$ of
Equation~\eqref{eq:G-hat}, and $\hat{S}(\bm k) = w\,\norm{\hat{\mat G}(\bm k)}_2^2$
by Equation~\eqref{eq:factored}. For $\bm k \ne \bm 0$ and each Gauss point
$g = 1, \dots, 4$, the $x$- and $y$-entries of $\hat{\mat V}(\bm k)$, in local
states $|2g-2\rangle$ and $|2g-1\rangle$, are
\begin{equation}
\label{eq:V-hat-entries}
\begin{aligned}
\hat V_{2g-2}(\bm k) &= \mathrm{i}\,e^{\mathrm{i}(\phi_1+\phi_2)}\,
\frac{\sin\phi_1\,(\cos\phi_2 + \mathrm{i}\,\eta_g\sin\phi_2)}{\hat{S}(\bm k)^{1/2}},
\\
\hat V_{2g-1}(\bm k) &= \mathrm{i}\,e^{\mathrm{i}(\phi_1+\phi_2)}\,
\frac{\sin\phi_2\,(\cos\phi_1 + \mathrm{i}\,\xi_g\sin\phi_1)}{\hat{S}(\bm k)^{1/2}},
\end{aligned}
\end{equation}
with
\begin{equation}
\label{eq:S-hat-explicit}
\begin{aligned}
\hat{S}(\bm k) = 4\Bigl(&\sin^2\phi_1\cos^2\phi_2 + \cos^2\phi_1\sin^2\phi_2
\\ &+ \tfrac23\,\sin^2\phi_1\sin^2\phi_2\Bigr).
\end{aligned}
\end{equation}
Summing the squared moduli of Equation~\eqref{eq:V-hat-entries} over the four
Gauss points, with $\xi_g^2 = \eta_g^2 = 1/3$, gives $1$, as it must.

By Equation~\eqref{eq:V-fourier}, the dense $8N^2 \times N^2$ matrix
$\mat{V}$ is replaced by $N^2$ unit vectors of length $8$, one per wavenumber,
each given in closed form. The circuit does not prepare these vectors one at a
time: a single wavenumber-controlled unitary prepares $\hat{\mat V}(\bm k)$
for all $\bm k$ in superposition, computing it from the bits of $\bm k$ at a
cost polylogarithmic in $N$ (Section~\ref{sec:cost}); a method
from a table of the $N^2$ vectors needs a circuit whose size grows as $N^2$. For a wavenumber along one axis, the vector reduces to the one-dimensional
phase: at $\bm k = (1, 0)$, the $y$-gradients vanish,
$\hat{S} = 4\sin^2\phi_1$, and
\begin{equation*}
\hat{\mat{V}}(1, 0) = \tfrac12\,\mathrm{i}\,e^{\mathrm{i}\pi/N}\,
(1, 0, 1, 0, 1, 0, 1, 0)^T ,
\end{equation*}
the phase of Equation~\eqref{eq:V-hat-1D} shared equally by the
$x$-gradients at the four Gauss points.
For the tensor-product two-point Gauss rule in $d$ dimensions, the vector has
length $d\,2^d$.

\subsection{Circuit for the Isometry}
\label{sec:circuit-V}

Equation~\eqref{eq:V-fourier} has the form of an FFT-based classical solver:
a transform to wavenumbers, one small operation per wavenumber, and a
transform back. We build the circuit $\mathcal{U}_V$ of
Figure~\ref{fig:circuit-V} with the same three stages. In this subsection, we
define $\mathcal{U}_V$ and show that it block-encodes $\mat{V}$ with
subnormalization $1$. The argument does not depend on how the middle stage is
synthesized; its synthesis, precision, and cost are given in
Section~\ref{sec:cost}.

\paragraph{Registers.}
The circuit must hold the input and the output of $\mat{V}$. The input is a
nodal vector $\bm u$, with one value at each of the $N^2$ nodes. The output
$\mat{V}\bm u$ is a gradient field, with eight values at each of the $N^2$
elements: two gradient components at each of four Gauss points. Each value is
therefore identified by a position, a node or an element, and by a local
index. At an element, the local index selects one of the $8$ values stored
there; at a node, it is $0$. The circuit uses one register for each. The \emph{position register} has
$m = \log_2 N^2$ qubits and holds the index $\bm j$ of a node or an element.
The \emph{local register} holds the local index $\ell$. A nodal vector is stored in local state $|0\rangle$, with its $N^2$ values in
the position register, and a gradient field in local states $|0\rangle$ to
$|7\rangle$, where $|2g-2\rangle$ holds the $x$-gradient and $|2g-1\rangle$
the $y$-gradient at Gauss point $g$. Between the two Fourier transforms, the
position register holds a wavenumber $\bm k$ instead of a position. There the
local rotation of the second stage below maps local state $|0\rangle$ to the
eight entries of $\hat{\mat V}(\bm k)$ in Equation~\eqref{eq:V-hat-entries},
stored in local states $|0\rangle$ to $|7\rangle$ with the same labels. At $\bm k = \bm 0$, $\hat{\mat V}(\bm 0) = \mat 0$, but
$\mathcal{Q}(\bm 0)$ is unitary and must map $|0\rangle$ to some state;
mapping it into $|0\rangle$ to $|7\rangle$ would add a gradient field that
$\mat{V}$ does not have. The circuit therefore needs local states outside the
gradient values. They are provided by a fourth qubit: the local register has four
qubits and $16$ states, the states $|8\rangle$ to $|15\rangle$ are not used
by $\mat{V}$, and they are used in Section~\ref{sec:be-M} to carry the
constant field through $\mathcal{U}_\Sigma$. In one
dimension, the local register is a single qubit: $|0\rangle$ holds the value
and $|1\rangle$ the constant field.

\paragraph{The three stages.}
\begin{enumerate}[leftmargin=1.4em,itemsep=2pt,topsep=2pt]
\item \textbf{Transform:} $\mathcal{F}$ on the position register maps the
nodal values to one amplitude per wavenumber $\bm k$.
\item \textbf{Local rotation:} for each wavenumber, a $16 \times 16$
unitary $\mathcal{Q}(\bm k)$ acts on the local register. For
$\bm k \ne \bm 0$, its first column is $\hat{\mat V}(\bm k)$ padded with
eight zeros, that is,
$\mathcal{Q}(\bm k)|0\rangle = \sum_{\ell=0}^{7}\hat V_\ell(\bm k)\,|\ell\rangle$
with the entries of Equation~\eqref{eq:V-hat-entries}, and it maps the states
$|0\rangle$ to $|7\rangle$ and the states $|8\rangle$ to $|15\rangle$ each
to themselves. Such a unitary exists because $\hat{\mat V}(\bm k)$ is a unit
vector (Section~\ref{sec:fourier-V}). For $\bm k = \bm 0$, where
$\hat{\mat V}(\bm 0) = \mat 0$, $\mathcal{Q}(\bm 0)$ maps $|0\rangle$ into
the unused states $|8\rangle$ to $|15\rangle$.
In Section~\ref{sec:cost}, $\mathcal{Q}(\bm k)$ is built with
fixed-point arithmetic on the bits of $\bm k$.
\item \textbf{Inverse transform:} $\mathcal{F}^{\dagger}$ on the position
register maps the wavenumbers back to elements.
\end{enumerate}
In one dimension, $\mathcal{Q}(k)$ multiplies $|0\rangle$ by the phase
$\mathrm{i}\,e^{\mathrm{i}\pi k/N}$ of Equation~\eqref{eq:V-hat-1D} for
$k \ne 0$ and maps $|0\rangle$ to $|1\rangle$ for $k = 0$.

\begin{figure}[H]
\centering
\resizebox{\linewidth}{!}{%
\begin{quantikz}
\lstick{local $|\ell\rangle$} & \qwbundle{4} & \qw & \qw
  & \gate[2]{\textstyle\sum_{\bm k}|\bm k\rangle\langle\bm k|
     \otimes \mathcal{Q}(\bm k)} & \qw & \qw & \qw\\
\lstick{position $|\bm j\rangle$} & \qwbundle{m} & \gate{\mathcal{F}}
  & \midstick[1,brackets=none]{$|\bm k\rangle$}
  & & \midstick[1,brackets=none]{$|\bm k\rangle$} & \gate{\mathcal{F}^{\dagger}} & \qw
\end{quantikz}}
\caption{The transform $\mathcal{F}$ maps mesh nodes to wavenumbers, the
wavenumber-controlled unitary $\mathcal{Q}(\bm k)$ maps the input local
state $|0\rangle$ to $\hat{\mathbf{V}}(\bm k)$, and $\mathcal{F}^{\dagger}$
maps back to mesh nodes. The output is read in local states $|0\rangle$ to
$|7\rangle$. As in Qiskit, the least significant register, here the local
register, is drawn at the top.}
\label{fig:circuit-V}
\end{figure}

\paragraph{Block encoding.}
The three stages compose to the unitary
\begin{equation}
\label{eq:UV}
\mathcal{U}_V = (\mathcal{F}^{\dagger} \kron I_{16})
  \Bigl(\sum_{\bm k} |\bm k\rangle\langle\bm k| \kron \mathcal{Q}(\bm k)\Bigr)
  (\mathcal{F} \kron I_{16}) .
\end{equation}
The transforms act on the position register only, and for every $\bm k$ the
components of $\mathcal{Q}(\bm k)|0\rangle$ in local states $|0\rangle$ to
$|7\rangle$ are the entries of $\hat{\mat V}(\bm k)$, including
$\bm k = \bm 0$, where $\mathcal{Q}(\bm 0)|0\rangle$ lies in the states
$|8\rangle$ to $|15\rangle$ and $\hat{\mat V}(\bm 0) = \mat 0$. In words: if
a nodal vector $\bm u$ is placed in local state $|0\rangle$ and $\mathcal{U}_V$
is applied, the part of the output in local states $|0\rangle$ to $|7\rangle$ is
the gradient field $\mat{V}\bm u$. Written with the projectors onto these local
states, Equation~\eqref{eq:V-fourier} gives
\begin{equation}
\label{eq:proj-encoding}
\begin{gathered}
\mat{V} = \bigl(I \kron \mat{\Pi}_{\mathrm{out}}\bigr)\,\mathcal{U}_V\,
\bigl(I \kron \mat{\Pi}_{\mathrm{in}}\bigr),
\\
\mat{\Pi}_{\mathrm{in}} = |0\rangle\langle 0| ,
\quad
\mat{\Pi}_{\mathrm{out}} = \sum_{\ell=0}^{7} |\ell\rangle\langle \ell| ,
\end{gathered}
\end{equation}
for every mesh size $N$: $\mathcal{U}_V$ block-encodes $\mat{V}$ with
subnormalization $1$. The constant field leaves the circuit in the local
states $|8\rangle$ to $|15\rangle$, outside the block, which is used in Section~\ref{sec:be-M}. One flag
qubit converts the projectors to the ancilla form of Equation~\eqref{eq:be}
\citep{gilyen2019quantum}, so the product rule of
Section~\ref{sec:quantum-preconditioning} applies to the products formed
below. Only the first column of $\mathcal{Q}(\bm k)$ and the two invariant
sets of local states are used in the argument. The gate-level circuit of
Section~\ref{sec:cost} computes that column with fixed-point
arithmetic, to within an encoding error $\delta$ defined there.

\paragraph{Relation to fast inversion.}
The circuit $\mathcal{U}_V$ has the structure of the fast-inversion circuit
of Tong et al.~\cite[Eq.~(10)]{tong2021fast}: a unitary transform, an
operation at each eigenvalue that computes a value by reversible arithmetic
and uses it to control a rotation, and the inverse transform. The difference
is the following. Tong et al.\ place an amplitude proportional to
$1/D_{ii}$ on a flag qubit, so their subnormalization is at least the norm of
the encoded matrix; built the same way, $(\mat{K}^{+})^{1/2}$ would have a
subnormalization that grows as $N$, which is the separately encoded
preconditioner of Equation~\eqref{eq:product-barrier}. Here $\hat S(\bm k)^{-1/2}$ is
folded into the unit vector $\hat{\mat V}(\bm k)$, so $\mathcal{Q}(\bm k)$ is
unitary and the subnormalization of $\mat{V}$ is $1$.

\subsection{Block Encoding of the Preconditioned Operator}
\label{sec:be-M}

With $\mathcal{U}_V$ in hand, the remaining factor of
Equation~\eqref{eq:M-VSV} is $\tilde{\mat{\Sigma}}$, and it is the only place
where the microstructure enters. We encode it with the material reflection
$R_\chi = I - 2\,\mathrm{diag}(\chi)$ of Section~\ref{sec:block-encoding},
applied to the element index and extended by the identity over the local
register. The reflection $-R_\chi$ is $+1$ at the quadrature points of the
inclusion and $-1$ at those of the matrix, while $\tilde{\mat{\Sigma}}$ is $1$
in the inclusion and $1/\rho$ in the matrix. Hence
\begin{equation}
\label{eq:E-lcu}
\begin{gathered}
\tilde{\mat{\Sigma}} = \beta_0\, I + \beta_1\,(-R_\chi),
\\
\beta_0 = \tfrac{1}{2}\bigl(1 + 1/\rho\bigr),
\quad
\beta_1 = \tfrac{1}{2}\bigl(1 - 1/\rho\bigr),
\end{gathered}
\end{equation}
since $\beta_0 + \beta_1 = 1$ and $\beta_0 - \beta_1 = 1/\rho$. This is a
linear combination of two unitaries with nonnegative coefficients that sum
to $1$. The standard circuit for such a combination
\citep{gilyen2019quantum} prepares one ancilla in
$\sqrt{\beta_0}\,|0\rangle + \sqrt{\beta_1}\,|1\rangle$, applies $-R_\chi$
controlled on $|1\rangle$, and undoes the preparation. The result is a block
encoding $\mathcal{U}_\Sigma$ of $\tilde{\mat{\Sigma}}$ with subnormalization
$\beta_0 + \beta_1 = 1$.

We now combine the two encodings into
$\mathcal{U}_M = \mathcal{U}_V^{\dagger}\,\mathcal{U}_\Sigma\,\mathcal{U}_V$,
which uses two queries to $\mathcal{U}_V$ and one controlled application of
$R_\chi$. To find the matrix it block-encodes, apply it to a nodal vector $\bm u$ in local state
$|0\rangle$, and split $\bm u$ into its constant part and the rest.
\begin{enumerate}[leftmargin=1.4em,itemsep=2pt,topsep=2pt]
\item \textbf{The part orthogonal to the constant fields:} $\mathcal{U}_V$
maps it to the gradient field $\mat{V}\bm u$ in the local states $|0\rangle$
to $|7\rangle$, $\mathcal{U}_\Sigma$ multiplies each gradient value by the
scaled conductivity of its element, and $\mathcal{U}_V^{\dagger}$ maps the
result back to local state $|0\rangle$. The block of these three steps is
$\mat{V}^T\tilde{\mat{\Sigma}}\,\mat{V} = \mat{M}$
(Equation~\eqref{eq:M-VSV}), with subnormalization $1 \cdot 1 \cdot 1 = 1$.
\item \textbf{The constant part:} $\mathcal{U}_V$ routes it into the unused local
states $|8\rangle$ to $|15\rangle$, $\mathcal{U}_\Sigma$ multiplies it element by element by
$\tilde{\mat{\Sigma}}$, and $\mathcal{U}_V^{\dagger}$ returns the constant
component of the product to local state $|0\rangle$. That component is the
constant field times $\bar{\sigma} = v_f + (1-v_f)/\rho$, the average of
$\tilde{\mat{\Sigma}}$ over the cell; for the square inclusion at $\rho = 10$,
$\bar{\sigma} = 0.325$.
\end{enumerate}
Therefore, the matrix that $\mathcal{U}_M$ block encodes is
\begin{equation}
\label{eq:M-tilde}
\tilde{\mat{M}} = \mat{M} + \bar{\sigma}\,\mat{\Pi}_0 ,
\end{equation}
with $\mat{\Pi}_0$ the projector onto the constant fields. 

Since $1/\rho \le \bar{\sigma} \le 1$, every eigenvalue of $\tilde{\mat{M}}$ lies
in $[1/\rho, 1]$, and the condition number seen by QSVT
(Equation~\eqref{eq:kappa-circuit}) is
\begin{equation}
\label{eq:kappa-M}
\kappa \;=\; 1\cdot\norm{\tilde{\mat{M}}^{-1}}_2 \;\le\; \rho
\end{equation}
for every mesh size $N$ and every phase arrangement $\chi$.

\subsection{Block Encoding of the Inverse}
\label{sec:inverse}
 QSVT applied to $\mathcal{U}_M$ therefore needs a polynomial that approximates $1/x$ only
on $1/\rho \le |x| \le 1$, and we use the odd polynomial $q_t$ of
Section~\ref{sec:qsvt} at $\kappa = \rho$. Gribling et
al.~\cite{gribling2021improving} bound the error of $q_t$ with a
Chebyshev estimate, and the same estimate gives the residual
$|x\,q_t(x) - 1| \le 2e^{-(d_t+1)/\rho}$ on that interval, where $d_t$ is the
degree, so a residual tolerance $\epsilon_p$ is met with
\begin{equation}
\label{eq:degree}
d_t \le \rho \ln(2/\epsilon_p) + 1
\end{equation}
queries to $\mathcal{U}_M$ \citep{gilyen2019quantum}. The subnormalization
$C = \max_{-1 \le x \le 1}|q_t(x)|$ is proportional to $\rho$ and grows as
$\sqrt{\ln(1/\epsilon_p)}$; over $1.5 \le \rho \le 10^{4}$ and
$10^{-2} \ge \epsilon_p \ge 10^{-12}$ we measure
$C \le \rho\sqrt{\tfrac12\ln(2/\epsilon_p)}$, an observation, not a theorem.

There are two sources of approximation errors: $q_t$ matches $1/x$
only to the tolerance $\epsilon_p$, and  $\mathcal{U}_V$
differs from an exact one by the encoding error $\delta$ of
Section~\ref{sec:cost}. Each of the $d_t$ queries applies
$\mathcal{U}_V$ and $\mathcal{U}_V^{\dagger}$ once, so the circuit
block-encodes a matrix within $2d_tC\delta$. For a heat source $\bm f$ orthogonal to the constant
fields, the computed solution
$\bm u_Q$ differs from the exact $\bm u = (\Kchi)^{+}\bm f$ by
\begin{equation}
\label{eq:energy-error}
\frac{\norm{\bm u_Q - \bm u}_{\Kchi}}{\norm{\bm u}_{\Kchi}}
\le \epsilon_p + 2d_tC\delta .
\end{equation}

\subsection{Gate-Level Circuit, Precision, and Cost}
\label{sec:cost}

The construction so far holds for any unitary $\mathcal{Q}(\bm k)$ whose
first column is $\hat{\mat V}(\bm k)$ (see Section \ref{sec:circuit-V}). A circuit must compute that column
from the bits of $\bm k$, coherently for all wavenumbers, and this
computation, not the Fourier transforms, sets the cost of the proposed method. In this
subsection, we separate the wavenumber dependence of $\hat{\mat V}(\bm k)$ into a real
vector of length $4$,
state the circuit that prepares it, which is constructed in
Appendix~\ref{app:multiplexed}, and then choose its precision and state
its cost.

\paragraph{Separating the wavenumber dependence.} The circuit does not have to
prepare $\hat{\mat V}(\bm k)$ directly. We factor it into three pieces: a
scalar phase, a fixed unitary independent of $\bm k$, and a real unit vector
of length $4$ that carries the whole dependence on $\bm k$. Only the third
piece needs arithmetic on the bits of $\bm k$; the first is $2n$ phase gates
and the second is one fixed circuit. With the abbreviations
\begin{equation}
\label{eq:abz}
a = \sin\phi_1\cos\phi_2,
\quad
b = \cos\phi_1\sin\phi_2,
\quad
\zeta = \sin\phi_1\sin\phi_2,
\end{equation}
the entries of Equation~\eqref{eq:V-hat-entries} at Gauss point $g$ are
$a + \mathrm{i}\,\eta_g\zeta$ and $b + \mathrm{i}\,\xi_g\zeta$, times the
common factor $\mathrm{i}\,e^{\mathrm{i}(\phi_1+\phi_2)}/\hat S(\bm k)^{1/2}$,
and Equation~\eqref{eq:S-hat-explicit} reads
$\hat S(\bm k) = 4(a^2 + b^2 + \tfrac23\zeta^2)$. Each entry is the sum of a
part that is the same at the four Gauss points, $a$ or $b$, and a part
proportional to $\eta_g$ or $\xi_g$. Only four patterns therefore occur in
$\hat{\mat V}(\bm k)$, whatever $\bm k$ is: the mean $x$- and $y$-gradients
and the two hourglass patterns. Written in the local states $|0\rangle$ to
$|7\rangle$, they are the columns of
\begin{equation}
\label{eq:t-patterns}
\mat{T}_4 = \begin{bmatrix}
\tfrac{\mathrm{i}}{2} & 0 & -\tfrac{\sqrt3}{2}\,\eta_g & 0\\[3pt]
0 & \tfrac{\mathrm{i}}{2} & 0 & -\tfrac{\sqrt3}{2}\,\xi_g
\end{bmatrix}_{g = 1,\dots,4},
\end{equation}
where the subscript $g$ denotes the four $2 \times 4$ blocks, one per Gauss
point, stacked into an $8 \times 4$ matrix, as in
Equation~\eqref{eq:G-hat}. Its columns are orthonormal, since
$\xi_g^2 = \eta_g^2 = 1/3$ and $\sum_g \xi_g = \sum_g \eta_g = 0$, so they
span a $4$-dimensional subspace of the $8$-dimensional gradient block, and
$\hat{\mat V}(\bm k)$ lies in that subspace for every $\bm k$. They are the
first four columns of a fixed unitary on the local register,
\begin{equation}
\label{eq:T-matrix}
\mat{T} = \begin{bmatrix}
\mat{T}_4 & \mat{T}_{\perp} & \mat{0}\\[2pt]
\mat{0} & \mat{0} & I_8
\end{bmatrix},
\end{equation}
where the columns of the $8 \times 4$ block $\mat{T}_{\perp}$ are any
orthonormal completion of the columns of $\mat{T}_4$ to a basis of the
gradient block, and $I_8$ is the identity on the unused states $|8\rangle$ to
$|15\rangle$. The completion is arbitrary: $\mat{V}$ depends on
$\mathcal{Q}(\bm k)$ only through its first column
(Section~\ref{sec:circuit-V}). Substituting Equation~\eqref{eq:abz} gives
\begin{equation}
\label{eq:Vhat-c}
\begin{gathered}
\hat{\mat V}(\bm k) = e^{\mathrm{i}(\phi_1+\phi_2)}\,\mat{T}\,
\frac{\bm c(\bm k)}{\norm{\bm c(\bm k)}_2},
\\
\bm c = \Bigl(a,\ b,\ \frac{\zeta}{\sqrt3},\ \frac{\zeta}{\sqrt3}\Bigr),
\quad
\norm{\bm c(\bm k)}_2 = \tfrac12\,\hat S(\bm k)^{1/2},
\end{gathered}
\end{equation}
so $\bm c(\bm k)$ holds the coordinates of $\hat{\mat V}(\bm k)$ in the basis
of Equation~\eqref{eq:t-patterns}, up to the phase. Its last two entries are
equal, the two hourglass patterns carrying the same weight, so three numbers
determine $\bm c(\bm k)$, and two after normalization. The phase is $2n$
phase gates on the bits of $\bm k$, with $n = \log_2 N$, and $\mat{T}$ is a
fixed circuit on three qubits of the local register.

\paragraph{The wavenumber-controlled preparation.} What remains is a
different real unit vector $\bm c(\bm k)/\norm{\bm c(\bm k)}_2$ for each of
the $N^2$ wavenumbers, a \emph{multiplexed state preparation}.
Appendix~\ref{app:multiplexed} constructs it in \emph{fixed-point
arithmetic}, each real number held in a register as an integer $Z$ and read
as $Z/2^F$, so that the number $F$ of \emph{fractional bits} sets the
resolution $2^{-F}$ of every computed value, and compares it with
standard state preparation multiplexed over $\bm
k$~\citep{mottonen2004transformation} and with a table lookup in a quantum
read-only memory~\citep{babbush2018encoding,low2024trading,berry2019qubitization}.
Both alternatives need $O(N^2)$ classical operations and a circuit whose
description grows as $N^2$ or $N$; the arithmetic is the only method whose
cost depends on $N$ only through $F$, and it is the one used here.

We use the main result of Appendix~\ref{app:multiplexed},
Proposition~\ref{prop:mx-main}: for $N = 2^n \ge 4$
and $F \ge n$ fractional bits, there is a circuit that prepares
$\bm c(\bm k)/\norm{\bm c(\bm k)}_2$ for all wavenumbers at once, with a
number of Toffoli gates that grows as $F^2$ and $2F + 3$ single-qubit
rotations, independent of $N$ except through $F$, on a work register of
$2F^2 + 35F + 47$ qubits that returns to $|0\rangle$ at every wavenumber, and
with an error $\delta(\bm k)$ of at most a constant times $F N\,2^{-F}$
(Proposition~\ref{prop:mx-error}). With the phase
$e^{\mathrm{i}(\phi_1+\phi_2)}$ and the fixed unitary $\mat{T}$ of
Equation~\eqref{eq:Vhat-c}, this circuit gives the gate-level
$\tilde{\mathcal{Q}}(\bm k)$. At $\bm k = \bm 0$, one $X$ gate on the most
significant local qubit, controlled on $\bm k = \bm 0$ at a cost of
$2(2n - 1)$ Toffoli gates, sends $|0\rangle$ into the unused states
$|8\rangle$ to $|15\rangle$; for $\bm k \ne \bm 0$,
$\tilde{\mathcal{Q}}(\bm k)$ keeps the states $|0\rangle$ to $|7\rangle$ and
$|8\rangle$ to $|15\rangle$ apart, as Section~\ref{sec:circuit-V} requires.
Since $\mat{T}$ is unitary, $\tilde{\mathcal{Q}}(\bm k)|0\rangle$ differs
from $\hat{\mat V}(\bm k)$ by the same $\delta(\bm k)$. The two are unit
vectors, so a rotation in their plane maps one onto the other and differs
from the identity by $\delta(\bm k)$; composing it with
$\tilde{\mathcal{Q}}(\bm k)$ gives an exact $\mathcal{Q}(\bm k)$ of
Section~\ref{sec:circuit-V}. The gate-level $\mathcal{U}_V$ therefore lies
within $\delta = \max_{\bm k}\delta(\bm k)$, in the spectral norm, of a
circuit that block-encodes $\mat{V}$ exactly. We call $\delta$ the
\emph{encoding error} of $\mathcal{U}_V$.

\paragraph{Precision.} The encoding error is linear in $N$ because the vector
$\bm c(\bm k)$ is short near $\bm k = \bm 0$:
$\norm{\bm c(\bm k)}_2 \ge \sin(\pi/N)$ for $\bm k \ne \bm 0$, the arithmetic
computes the entries of $\bm c(\bm k)$ to an absolute error of order
$2^{-F}$, and the normalization divides that error by $\norm{\bm c(\bm k)}_2$
(Appendix~\ref{app:mx-error}). We choose $F$ from the measured bound
\begin{equation}
\label{eq:delta-bound}
\delta \le 6N\,2^{-F},
\end{equation}
an observation over $42$ simulated configurations with $N \le 256$
(Appendix~\ref{app:mx-error}), which we extrapolate to larger $N$.
By Equation~\eqref{eq:energy-error}, the error
of the inverse stays within $2\epsilon_p$ if $2d_tC\delta \le \epsilon_p$, and
requiring $6N\,2^{-F} \le \epsilon_p/(2d_tC)$ gives
\begin{equation}
\label{eq:F-choice}
F = \Bigl\lceil \log_2 \frac{12\,N\,d_t\,C}{\epsilon_p} \Bigr\rceil .
\end{equation}
For the square at $N = 1024$, $\rho = 10$, and $\epsilon_p = 10^{-6}$, this
gives $F = 45$: the mesh, the polynomial, and the tolerance contribute about
$10$, $11.3$, and $19.9$ bits, and the constant $12$ another $3.6$. At
$\epsilon_p = 10^{-3}$, it gives $F = 34$. Example \ref{ex:precision} illustrates the meaning of $F$.

\begin{example}[The meaning of $F$]
\label{ex:precision}
Let $N = 16$ and $F = 12$, so that the resolution of every register is
$\Delta = 2^{-12} = 2.4 \times 10^{-4}$, and take $\bm k = (3, 1)$, for which
$\phi_1 = 3\pi/16$ and $\phi_2 = \pi/16$. The exact entries of
Equation~\eqref{eq:abz} are $a = 0.54490$, $b = 0.16221$, and
$\zeta = 0.10839$; the arithmetic of Appendix~\ref{app:mx-part1} computes
them as $0.54419$, $0.16187$, and $0.10791$, in error by $2.9\Delta$,
$1.4\Delta$, and $2.0\Delta$. The resolution of the
registers therefore sets the absolute error of the entries, to a few units
of $\Delta$. Here $\norm{\bm c(\bm k)}_2 = 0.57537$, and
\begin{equation}
\label{eq:ex-chat}
\frac{\bm c(\bm k)}{\norm{\bm c(\bm k)}_2} =
(0.94703,\ 0.28192,\ 0.10876,\ 0.10876)
\end{equation}
against the prepared
$(0.94709,\ 0.28202,\ 0.10835,\ 0.10835)$, so that
$\delta(\bm k) = 5.9 \times 10^{-4}$, or $2.4\Delta$. The bound of Equation~\eqref{eq:delta-bound} is $6N\,2^{-F} = 0.023$ at these settings, about $40$ times larger. At the settings of this paper $F = 45$:
$\Delta = 2.8 \times 10^{-14}$, $\delta \le 6 \cdot 1024 \cdot 2^{-45} =
1.7 \times 10^{-10}$, and $2d_tC\delta = 8.9 \times 10^{-7}$, just inside
$\epsilon_p = 10^{-6}$ as Equation~\eqref{eq:F-choice} requires.
\end{example}

\paragraph{Cost.} With $F$ fixed, we count the gates that are expensive on a
fault-tolerant quantum computer: Toffoli gates, and rotations of a single
qubit by an arbitrary angle, which we count as they appear in the circuit,
before they are compiled into Clifford and T gates. The Fourier transforms
use no Toffoli gates, so all of them are in the arithmetic of
Appendix~\ref{app:multiplexed}, $96$ percent in its four CORDIC runs, and their number grows as $F^2$,
that is, as $(\log N + \log(1/\delta))^2$. For the square at $N = 1024$ and
$F = 45$:
\begin{itemize}[leftmargin=1.4em,itemsep=2pt,topsep=2pt]
\item \textbf{One application of $\mathcal{U}_V$:} $44{,}462$ Toffoli gates on
$5{,}696$ qubits: the $20$ qubits of the position register, the $4$ of the
local register, and a work register of $5{,}672$ qubits, which holds the
intermediate values of the arithmetic and returns to $|0\rangle$.
\item \textbf{One query to $\mathcal{U}_M$:} $88{,}926$ Toffoli gates and $588$
rotations on $5{,}699$ qubits. A query applies $\mathcal{U}_V$ twice and the
controlled material reflection once, and adds three single ancillas, for
$\mathcal{U}_\Sigma$, the flag of Section~\ref{sec:circuit-V}, and QSVT. Its
rotations are, for each $\mathcal{U}_V$, the $2F + 3 = 93$ of the
arithmetic, the $2n = 20$ phase gates, and the $180$
controlled phase gates of the two Fourier transforms, which apply the
one-dimensional quantum Fourier transform to each coordinate register with
$n(n-1)/2$ controlled phases each~\citep{suresh2026appliedqc}; and the $2$
rotations of $\mathcal{U}_\Sigma$.
\end{itemize}
One application of the inverse therefore costs
$O\bigl(\rho\log(1/\epsilon_p)\,(\log N + \log(1/\delta))^2\bigr)$ Toffoli
gates, the $d_t$ queries of Equation~\eqref{eq:degree} at the cost per query
above; Section~\ref{sec:exp-resources} gives the count at $N = 1024$.

\section{Implementation}
\label{sec:implementation}

The constructions of this paper are implemented in \pkgname{}, an
open-source repository available at \url{\repourl}, which also holds the code
of Appendix~\ref{app:multiplexed}. It requires NumPy, SciPy,
Matplotlib for the figures, and Qiskit~\citep{javadi2024quantum} and pyqsp for
the circuit of Listing~\ref{lst:qiskit}. The repository has three layers.
\begin{itemize}[leftmargin=1.4em,itemsep=2pt,topsep=2pt]
\item \textbf{The package \texttt{\libname}:} \texttt{fourier\_symbol.py}
holds the closed-form eigenvalues of $\mat{K}$, functions of $\mat{K}$ by the
FFT, and $\mat{K}$ assembled element by element; \texttt{isometry.py} holds
$\hat{\mat V}(\bm k)$ in closed form, $\mat{V}$ assembled element by element,
the unitary $\mat{T}$ of Equation~\eqref{eq:Vhat-c}, and the gate-level
circuit $\mathcal{U}_V$ of Section~\ref{sec:cost}, whose arithmetic, the
circuit of Appendix~\ref{app:multiplexed}, is in the subpackage
\texttt{arithmetic}: \texttt{revcirc.py} holds the reversible fixed-point
arithmetic and a simulator that runs a circuit bit by bit on every basis input
at once, \texttt{symbol.py} and \texttt{tree.py} hold Parts 1 and 2 of the
circuit, and \texttt{alternatives.py} holds the cost estimates of the
multiplexor and the table lookup.
\item \textbf{The test cells and the inverse:} \texttt{cells.py} holds the
three cells of Section~\ref{sec:exp-cells}, their material oracles, and the
operators applied without dense matrices; \texttt{inversion.py} holds the
inverse, the choice of $F$, and the resource counts.
\item \textbf{Scripts:} in the folder \texttt{papers/inversion},
\texttt{makeTables.py} and \texttt{makeFigures.py} write every table and
figure of the paper, Appendix~\ref{app:multiplexed} included, after asserting
every number the text quotes about them.
\end{itemize}

The simplest use of the package is on the square inclusion of the preceding
sections. In Listing~\ref{lst:inverse}, we build the square cell at $N = 256$ and
$\rho = 10$, compute the extreme eigenvalues of $\tilde{\mat{M}}$, apply the
inverse of Section~\ref{sec:inverse} at $\epsilon_p = 10^{-6}$ to the source of
the cell problem for the macroscopic gradient $(1, 0)$, and then compute the
resources of the inverse at $N = 1024$. Each quantity is a single
function call in the package.

\begin{lstlisting}[caption={The inverse on the square at $N = 256$, and its resources at $N = 1024$.},label={lst:inverse},columns=fixed]
from quantumpoissonhomogenization import (Cell, extreme_eigs, inverse_error,
                                          inverse_resources)

cell = Cell("square", N=256, rho=10.0)             # 65,536 unknowns
lo, hi = extreme_eigs(cell)                        # extreme eigenvalues of Mtilde
err, d = inverse_error(cell, eps_p=1e-6)           # q_t(Mtilde), matrix-free
print(f"kappa(Mtilde)    = {hi / lo:.6f}")
print(f"degree d_t       = {d}")
print(f"energy error     = {err:.2e}")

v = inverse_resources("square", N=1024, rho=10.0, eps_p=1e-6)
print(f"F                = {v['F']}")
print(f"Toffoli / query  = {v['toffoli_per_query']:,}")
print(f"Toffoli, inverse = {v['toffoli_total']:.2e}")
print(f"qubits           = {v['qubits']:,}")
\end{lstlisting}

\begin{lstlisting}[numbers=none,frame=none,backgroundcolor=\color{white},columns=fixed,language={}]
----- output --------------
kappa(Mtilde)    = 10.000000
degree d_t       = 145
energy error     = 5.03e-07
F                = 45
Toffoli / query  = 88,926
Toffoli, inverse = 1.29e+07
qubits           = 5,699
\end{lstlisting}
\noindent The condition number equals $\rho$, and the inverse at $d_t = 145$
meets $\epsilon_p = 10^{-6}$. These computations are repeated in
Section~\ref{sec:numerics} for three cells and meshes up to $N = 1024$
(Figure~\ref{fig:mesh}), and the resources are those of
Table~\ref{tab:resources}. With \texttt{"cross"} or \texttt{"laminate"} in
place of \texttt{"square"}, the same computations are run on the other two
cells.

In Listing~\ref{lst:inverse}, $q_t(\tilde{\mat{M}})$ is applied without
gates. In Listing~\ref{lst:uv}, we check the gate-level circuit itself: we
build $\mathcal{U}_V$ at $N = 16$ and $F = 24$, simulate it gate by gate for
all $256$ wavenumbers, and compare the vector that $\mathcal{Q}(\bm k)$
prepares with $\hat{\mat V}(\bm k)$ of
Equation~\eqref{eq:V-fourier}.

\begin{lstlisting}[caption={The gate-level $\mathcal{U}_V$ at $N = 16$ and $F = 24$, checked against $\hat{\mathbf{V}}(\boldsymbol{k})$ at every wavenumber.},label={lst:uv},columns=fixed]
import numpy as np
from quantumpoissonhomogenization import IsometryCircuit, vhat_exact

circ = IsometryCircuit(n=4, F=24)          # N = 16, F = 24 fractional bits
Q, checks = circ.simulate()                # Q(k)|0>, gate by gate, every k
N = circ.N
err = max(np.linalg.norm(Q[k1 * N + k2, :8, 0] - vhat_exact(N, k1, k2))
          for k1 in range(N) for k2 in range(N))
r = circ.resources()

print(f"wavenumbers      = {N * N}")
print(f"work register    = {checks['work_bits_left_set']} bits left set")
print(f"inputs restored  = {checks['inputs_restored']}")
print(f"max |error|      = {err:.1e}")
print(f"Toffoli gates    = {r['toffoli_total']:,}")
print(f"qubits           = {r['wires']:,}")
\end{lstlisting}

\begin{lstlisting}[numbers=none,frame=none,backgroundcolor=\color{white},columns=fixed,language={}]
----- output --------------
wavenumbers      = 256
work register    = 0 bits left set
inputs restored  = True
max |error|      = 1.1e-06
Toffoli gates    = 13,349
qubits           = 2,051
\end{lstlisting}
\noindent The work register returns to $|0\rangle$ at every wavenumber, and
the prepared vectors agree with $\hat{\mat V}(\bm k)$ to
$1.1 \times 10^{-6}$, below the bound $6N\,2^{-F} = 5.7 \times 10^{-6}$ of
Equation~\eqref{eq:delta-bound}.

In Listings~\ref{lst:inverse} and~\ref{lst:uv}, the inverse is applied as a
matrix polynomial, and the gate-level circuit is checked one stage at a time.
In Listing~\ref{lst:qiskit}, we build the entire inverse, from the right-hand
side to the output state, as one Qiskit circuit for the square at $N = 4$,
$\rho = 10$, and $\epsilon_p = 10^{-3}$. We then simulate it by state vector
and compare the output state with the finite element solution. The circuit is
described in Section~\ref{sec:exp-qiskit}.

\begin{lstlisting}[caption={The entire inverse as one Qiskit circuit at $N = 4$, simulated by state vector.},label={lst:qiskit},columns=fixed]
from quantumpoissonhomogenization.inversion import qiskit_inverse

r = qiskit_inverse(N=4, rho=10.0, eps_p=1e-3)    # square cell, 16 unknowns
print(f"degree d_t       = {r['d']}")
print(f"qubits           = {r['qubits']}")
print(f"CX gates         = {r['inverse'][1]:,}")
print(f"block of U_M     = {r['err_M']:.1e}")
print(f"energy error     = {r['energy']:.2e}")
print(f"P(success)       = {r['p_success']:.2e}")
\end{lstlisting}

\begin{lstlisting}[numbers=none,frame=none,backgroundcolor=\color{white},columns=fixed,language={}]
----- output --------------
degree d_t       = 77
qubits           = 13
CX gates         = 25,938
block of U_M     = 2.9e-14
energy error     = 6.40e-04
P(success)       = 2.44e-03
\end{lstlisting}
\noindent The circuit uses $13$ qubits, its block of $\mathcal{U}_M$ equals
$\tilde{\mat{M}}$ to machine precision, and the energy-norm error of its
output is below $\epsilon_p = 10^{-3}$.

Listings~\ref{lst:inverse} and~\ref{lst:uv} run from the folder \texttt{papers/inversion} in
a few seconds, and Listing~\ref{lst:qiskit} in about $15$ seconds. The two
scripts regenerate every table and figure of the paper:
\texttt{python makeTables.py} in about $3$ minutes and
\texttt{python makeFigures.py} in about $10$ minutes.
These results are reported in Section~\ref{sec:numerics}.

\section{Numerical Experiments}
\label{sec:numerics}

We now check the claims of Section~\ref{sec:proposed} on three periodic
cells, for meshes up to
$N = 1024$, that is, $N_{\mathrm{dof}} \approx 1.05 \times 10^6$ unknowns.
We address four questions.
\begin{enumerate}[leftmargin=1.4em,itemsep=2pt,topsep=2pt]
\item \textbf{Geometry:} is the condition number of $\tilde{\mat{M}}$ equal to
$\rho$ for cells of different topology (Section~\ref{sec:exp-cells})?
\item \textbf{Mesh:} does it stay at $\rho$ under refinement, and does the
inverse keep its accuracy at a fixed degree (Section~\ref{sec:exp-mesh})?
\item \textbf{Gates:} does the gate-level circuit encode $\tilde{\mat{M}}$ and
its inverse to the accuracy that the encoding error predicts
(Sections~\ref{sec:exp-be} and~\ref{sec:exp-inverse})?
\item \textbf{End to end:} does the entire circuit, from the right-hand side to
the output state, deliver the finite element solution
(Section~\ref{sec:exp-qiskit})?
\end{enumerate}
The logical cost of the inverse on the finest mesh is then given in
Section~\ref{sec:exp-resources}.

Every table and figure of this section is written by the two scripts of
Section~\ref{sec:implementation}. No dense matrix of size $N_{\mathrm{dof}}$
is formed on the large meshes. The quantities are computed as follows.
\begin{itemize}[leftmargin=1.4em,itemsep=2pt,topsep=2pt]
\item \textbf{Operators:} functions of $\mat{K}$ are applied by the FFT, and
$\mat{W}^{1/2}\mat{G}$ as a sparse matrix.
\item \textbf{Inverse:} $q_t(\tilde{\mat{M}})$ is applied by the Clenshaw
recurrence on the Chebyshev coefficients of $q_t$, which uses $d_t$
applications of $\tilde{\mat{M}}$.
\item \textbf{Eigenvalues:} the extreme eigenvalues are computed by Lanczos
iteration.
\item \textbf{Reference solutions:} the two-phase problem is solved by
preconditioned conjugate gradients to a relative residual of $10^{-13}$.
\item \textbf{Gate counts:} these come from the gate-level circuit of
Section~\ref{sec:cost} and from the material oracles of
Section~\ref{sec:exp-cells}.
\item \textbf{Dense checks:} on meshes up to $N = 64$, the encoding and the
inverse are also compared with dense matrices and a sparse direct solve.
\item \textbf{Qiskit circuit:} at $N = 4$ and $8$, the entire inverse is built
from Qiskit gates and simulated by state vector.
\end{itemize}

\subsection{Test Cells}
\label{sec:exp-cells}

\begin{figure}[H]
\centering
\includegraphics[width=0.975\linewidth]{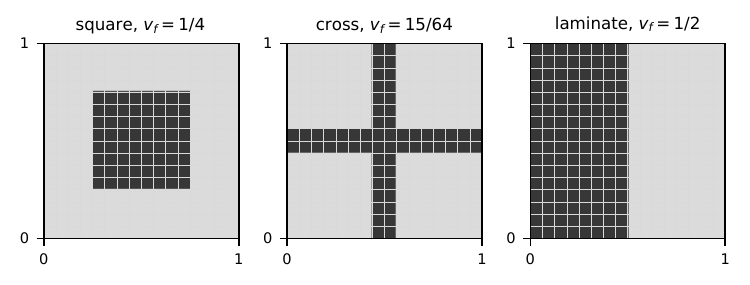}
\caption{The three test cells, drawn at $N = 16$; the inclusion phase is
dark. The square is an isolated conductor, the cross a connected one, and
the laminate an anisotropic cell.}
\label{fig:cells}
\end{figure}

We use the three cells of Figure~\ref{fig:cells}, chosen to differ in
topology. The square inclusion of the preceding sections is an isolated
conductor. The cross, of nearly the same volume fraction $v_f = 15/64$, is a
conductor connected across the cell in both directions. The laminate, with
layers normal to $x$, is anisotropic. All three cells are \emph{dyadic}: the
phase indicator $\chi$ is a Boolean function of the two (square), four
(cross), or one (laminate) leading bits of each coordinate. The controlled
material reflection of Section~\ref{sec:be-M} therefore costs $2$, $9$, and
$0$ Toffoli gates, independent of $N$. We check each of these circuits on
every element for $N = 16$ to $256$: it reproduces $2\chi - 1$ in every case.

For these three cells, Table~\ref{tab:geometry} reports the extreme
eigenvalues of $\tilde{\mat{M}}$ at $N = 256$, for $\rho = 10$ and $\rho = 10^4$. In all six cases,
$\lambda_{\min} = 1/\rho$ and $\lambda_{\max} = 1$ to six digits, so
$\kappa(\tilde{\mat{M}}) = \rho$. The bound of Equation~\eqref{eq:kappa-M} is
therefore attained for an isolated, a connected, and a layered conductor
alike. This is the behavior predicted in
Section~\ref{sec:reference-preconditioning}: a field whose gradient vanishes in every inclusion element gives the
eigenvalue $1/\rho$, and a field whose gradient vanishes in every matrix
element gives the eigenvalue $1$.

\begin{table*}[t]
\caption{Extreme eigenvalues of $\tilde{\mat{M}}$ for the three cells at
$N = 256$. ``$R_\chi$ Toffoli'' is the number of Toffoli gates of the controlled
material reflection; $\rho\,\lambda_{\min} = \lambda_{\max} = 1$ means that the
bound $\kappa(\tilde{\mat{M}}) \le \rho$ of Equation~\eqref{eq:kappa-M} is
attained.}
\label{tab:geometry}
\centering
\small
\setlength{\tabcolsep}{4pt}
\begin{tabular}{lrrrrrr}
\toprule
geometry & $v_f$ & $R_\chi$ Toffoli & $\rho$ & $\rho\,\lambda_{\min}$ & $\lambda_{\max}$ & $\kappa(\tilde{\mat{M}})/\rho$\\
\midrule
square & $0.25$ & $2$ & $10$ & $1.000000$ & $1.000000$ & $1.000000$\\
 &  &  & $10^{4}$ & $1.000000$ & $1.000000$ & $1.000000$\\
cross & $0.234$ & $9$ & $10$ & $1.000000$ & $1.000000$ & $1.000000$\\
 &  &  & $10^{4}$ & $1.000000$ & $1.000000$ & $1.000000$\\
laminate & $0.5$ & $0$ & $10$ & $1.000000$ & $1.000000$ & $1.000000$\\
 &  &  & $10^{4}$ & $1.000000$ & $1.000000$ & $1.000000$\\
\bottomrule
\end{tabular}

\end{table*}

\subsection{Mesh Independence}
\label{sec:exp-mesh}

\begin{figure*}[t]
\centering
\includegraphics[width=0.8\linewidth]{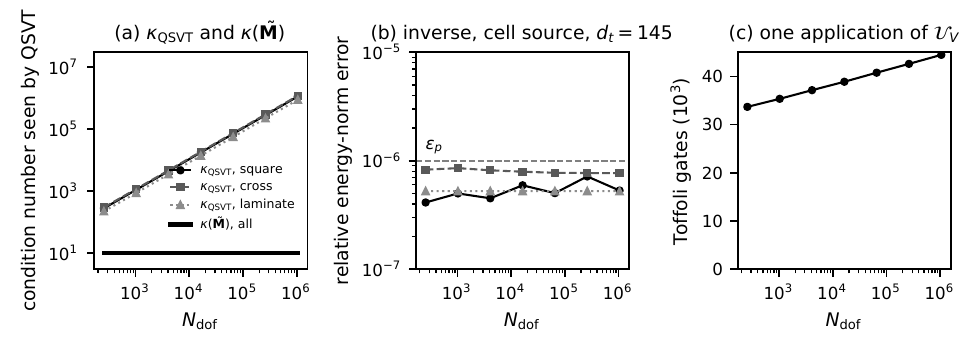}
\caption{Mesh refinement at $\rho = 10$. (a) Condition number seen by QSVT
for the two-phase operator ($\alpha = 16/3$) of the three cells, and
$\kappa(\tilde{\mat{M}})$, which is the same for all three cells. (b)
Relative energy-norm error of the inverse at the fixed degree $d_t = 145$
($\epsilon_p = 10^{-6}$), for the cell-problem source of each cell, with the
exact $\mathcal{Q}(\bm k)$. (c) Toffoli gates of one application of $\mathcal{U}_V$, with $F$ set by
Equation~\eqref{eq:F-choice}; rotations are counted in
Table~\ref{tab:resources}.}
\label{fig:mesh}
\end{figure*}

In Figure~\ref{fig:mesh}, the mesh is refined from $N = 16$ to $1024$ at
$\rho = 10$ to test three predictions of Sections~\ref{sec:background}
and~\ref{sec:proposed}, one per panel.

\paragraph{(a) Condition number.} For the two-phase operator, the condition
number that QSVT sees is $\alpha\,\norm{(\Kchi/\sigma_{\max})^{+}}_2$ with
$\alpha = 16/3$. It grows by a factor of $4.00$ per doubling of $N$ on the
finest meshes and reaches $1.1 \times 10^{6}$ for the square at $N = 1024$,
the $\Theta(N^2)$ growth of Section~\ref{sec:qsvt}. The condition number of
$\tilde{\mat{M}}$ equals $\rho = 10$ to ten digits for every $N$ and every
cell.

\paragraph{(b) Accuracy of the inverse.} Because $\kappa(\tilde{\mat{M}})$
does not change with $N$, we fix the degree at $d_t = 145$, the
degree of Section~\ref{sec:inverse} for $\epsilon_p = 10^{-6}$, and
apply the inverse to the cell-problem source of
Section~\ref{sec:exp-inverse}. The relative energy-norm error stays between $4.1 \times 10^{-7}$
and $8.6 \times 10^{-7}$ for all three cells and all meshes, below
$\epsilon_p$ as Equation~\eqref{eq:energy-error} requires. The error does not
grow with $N$, so the degree does not need to grow either.

\paragraph{(c) Cost of one application of $\mathcal{U}_V$.} What grows with
$N$ is the cost of $\mathcal{U}_V$. By Equation~\eqref{eq:F-choice}, $F$
grows by one bit per doubling of $N$, from $F = 39$ at $N = 16$ to $F = 45$ at
$N = 1024$, and the Toffoli count grows quadratically in $F$
(Section~\ref{sec:cost}), from $33{,}680$ to $44{,}462$. The cost
per query therefore grows with $N$ only through $F$, which grows as
$\log_2 N$, while the number of queries stays fixed.

\subsection{Gate-Level Encoding}
\label{sec:exp-be}

The block of
$\mathcal{U}_M = \mathcal{U}_V^{\dagger}\mathcal{U}_\Sigma\mathcal{U}_V$ was
derived in Section~\ref{sec:be-M} with the exact $\mathcal{Q}(\bm k)$. We now form that block from the simulated gate-level circuit,
\begin{equation*}
(I \kron \mat{\Pi}_{\mathrm{in}})\,\mathcal{U}_V^{\dagger}\,
\mathcal{U}_\Sigma\,\mathcal{U}_V\,(I \kron \mat{\Pi}_{\mathrm{in}}),
\end{equation*}
including the routing of $\bm k = \bm 0$ into the unused local states, and
compare it with $\tilde{\mat{M}}$ for the square at $\rho = 10$.
Table~\ref{tab:be-M} reports the spectral-norm difference for $N = 8$ and
$16$ and $F = 24$ and $32$. In every case the difference lies below $2\delta$,
as two applications of $\mathcal{U}_V$ allow, and it falls with $F$. The
smallest and largest eigenvalues of the block equal $1/\rho$ and $1$ to $13$
digits, and the constant field appears as the eigenvalue
$\bar{\sigma} = 0.325$, inside that interval.

\begin{table}[t]
\caption{Spectral-norm difference between the block of
$\mathcal{U}_V^{\dagger}\mathcal{U}_\Sigma\mathcal{U}_V$ on the input subspace,
computed from the gate-level circuit with $F$ fractional bits, and
$\tilde{\mat{M}}$, for the square at $\rho = 10$. The difference is the encoding
error of the circuit, not a discretization error; it falls with $F$ and lies
below $2\delta$, twice the encoding error of $\mathcal{U}_V$ at the same $N$
and $F$.}
\label{tab:be-M}
\centering
\begin{tabular}{rrrr}
\toprule
$N$ & $F$ & difference & $2\delta$\\
\midrule
$8$ & $24$ & $1.8 \times 10^{-7}$ & $2.2 \times 10^{-6}$\\
$8$ & $32$ & $6.7 \times 10^{-10}$ & $7.4 \times 10^{-9}$\\
$16$ & $24$ & $2.1 \times 10^{-7}$ & $2.2 \times 10^{-6}$\\
$16$ & $32$ & $8.0 \times 10^{-10}$ & $1.0 \times 10^{-8}$\\
\bottomrule
\end{tabular}

\end{table}

\subsection{The Inverse on Four Heat Sources}
\label{sec:exp-inverse}

We apply the inverse of Section~\ref{sec:inverse}, with $\rho = 10$ and
$\epsilon_p = 10^{-6}$, to four heat sources on the square cell, each with its
mean removed.
\begin{itemize}[leftmargin=1.4em,itemsep=2pt,topsep=2pt]
\item \textbf{Smooth:} the nodal source $\sin 2\pi x\cos 2\pi y$.
\item \textbf{Random:} a source with independent standard normal entries.
\item \textbf{Cell problem:} the source
$\bm f = -\mat{G}^T\mat{W}\mat{\Sigma}\,\bar{\bm e}$, with $\bar{\bm e}$ the
gradient $(1, 0)$ at every quadrature point. This is the source of the cell
problem of homogenization for the macroscopic gradient $(1, 0)$; it is
concentrated at the interface between the phases and orthogonal to the
constant fields, since $\mat{G}$ maps them to zero.
\item \textbf{Worst case:} $\mat{K}^{1/2}\bm v$, with $\bm v$ an eigenvector of
$\mat{M}$ at $1/\rho$. For this source the polynomial part of
Equation~\eqref{eq:energy-error} is attained.
\end{itemize}
For each source, the reference is the finite element solution from a sparse
direct solve. Table~\ref{tab:inverse} reports the relative error in the energy norm. It is
at most $\epsilon_p$ in every case and does not grow with $N$, and for the
worst-case source it equals the residual $8.69 \times 10^{-7}$ of $q_t$. We
then replace $\tilde{\mat{M}}$ by the block of the gate-level circuit at
$N = 8$ and $16$. The errors change by at most $1.2 \times 10^{-7}$ for
$F = 24$ and $2.7 \times 10^{-9}$ for $F = 32$, within the encoding term of
Equation~\eqref{eq:energy-error}.

\begin{table*}[t]
\caption{Relative energy-norm error of the inverse with $\rho = 10$ and
$\epsilon_p = 10^{-6}$ for four heat sources on the square cell, with the
degree $d_t$. The errors are those of the polynomial applied to the exact
$\tilde{\mat{M}}$; they do not include the encoding error of the circuit.}
\label{tab:inverse}
\centering
\small
\begin{tabular}{rrrrrr}
\toprule
$N$ & $d_t$ & smooth & random & cell & worst\\
\midrule
$8$ & $145$ & $8.4 \times 10^{-7}$ & $8.4 \times 10^{-7}$ & $4.7 \times 10^{-7}$ & $8.7 \times 10^{-7}$\\
$16$ & $145$ & $8.4 \times 10^{-7}$ & $8.6 \times 10^{-7}$ & $4.1 \times 10^{-7}$ & $8.7 \times 10^{-7}$\\
$32$ & $145$ & $8.5 \times 10^{-7}$ & $8.6 \times 10^{-7}$ & $5.0 \times 10^{-7}$ & $8.7 \times 10^{-7}$\\
$64$ & $145$ & $8.5 \times 10^{-7}$ & $8.6 \times 10^{-7}$ & $4.5 \times 10^{-7}$ & $8.7 \times 10^{-7}$\\
\bottomrule
\end{tabular}

\end{table*}
\subsection{The Inverse as One Qiskit Circuit}
\label{sec:exp-qiskit}

\begin{table*}[t]
\caption{Logical resources of one application of the inverse
$q_t(\tilde{\mat{M}})$ for the square at $N = 1024$ and $\rho = 10$: polynomial
tolerance $\epsilon_p$, degree $d_t$ (queries to $\mathcal{U}_M$), subnormalization
$C$, fractional bits $F$ from Equation~\eqref{eq:F-choice}, Toffoli gates of one
application of $\mathcal{U}_V$ and of one query, single-qubit rotations per
query, total Toffoli gates, and logical qubits. Rotations are counted before
synthesis into Clifford and T gates.}
\label{tab:resources}
\centering
\footnotesize
\setlength{\tabcolsep}{3pt}
\begin{tabular}{rrrrrrrrr}
\toprule
$\epsilon_p$ & $d_t$ & $C$ & $F$ & Toffoli/$\mathcal{U}_V$ & Toffoli/query & rotations/query & Toffoli & qubits\\
\midrule
$10^{-3}$ & $77$ & $13.2$ & $34$ & $25{,}925$ & $51{,}852$ & $544$ & $4.0 \times 10^{6}$ & $3{,}576$\\
$10^{-6}$ & $145$ & $17.6$ & $45$ & $44{,}462$ & $88{,}926$ & $588$ & $1.3 \times 10^{7}$ & $5{,}699$\\
\bottomrule
\end{tabular}

\end{table*}

In Sections~\ref{sec:exp-be} and~\ref{sec:exp-inverse}, the encoding and the
inverse are checked separately, and the inverse is applied as a matrix
polynomial. We now build the entire inverse, from the right-hand side to the
output state, as one circuit of Qiskit gates~\cite{javadi2024quantum} and
simulate it by state vector. The fixed-point arithmetic of
Section~\ref{sec:cost} cannot be part of such a simulation: its work register
alone has $2F^2 + 35F + 47 = 755$ qubits at $F = 12$. We therefore build
$\mathcal{Q}(\bm k)$ as the multiplexor of Appendix~\ref{app:multiplexed}, whose
angles are computed classically for all $N^2$ wavenumbers and which needs no
work register; the arithmetic is checked separately
(Listing~\ref{lst:uv} and Table~\ref{tab:be-M}). The settings are those of
Section~\ref{sec:inverse}: the square at $\rho = 10$ and
$\epsilon_p = 10^{-3}$, so that $d_t = 77$ and $C = 13.2$, on meshes with
$N = 4$ and $N = 8$. The construction is
Listing~\ref{lst:qiskit}, and the circuits are built in
\texttt{inversion.py} of the repository.

The same circuits prepare the right-hand side, with no classical
preprocessing. The cell-problem source of Section~\ref{sec:exp-inverse} is
$\bm f = -\mat{G}^T\mat{W}\mat{\Sigma}\,\bar{\bm e}$; splitting
$\mat{W} = \mat{W}^{1/2}\mat{W}^{1/2}$ and using Equation~\eqref{eq:V}, the
preconditioned right-hand side is
\begin{equation}
\label{eq:psi-rhs}
\begin{aligned}
\bm\psi = (\mat{K}^{+})^{1/2}\bm f
&= -\bigl[(\mat{K}^{+})^{1/2}\mat{G}^T\mat{W}^{1/2}\bigr]\,
\mat{W}^{1/2}\mat{\Sigma}\,\bar{\bm e}
\\
&= -\sqrt{w}\,\sigma_{\max}\,\mat{V}^T\tilde{\mat{\Sigma}}\,\bar{\bm e} ,
\end{aligned}
\end{equation}
so that Hadamard gates, one application of $\mathcal{U}_\Sigma$ and one of
$\mathcal{U}_V^{\dagger}$ deliver $\bm\psi$, flagged on two ancillas. QSVT
then applies $q_t$ to $\mathcal{U}_M$ with the phases of Dong et
al.~\cite{dong2021efficient}, and the whole circuit uses $13$ qubits at
$N = 4$ and $15$ at $N = 8$: $2n$ position qubits, $4$ local qubits, and $5$
ancillas.

Before transpiling, each stage is compared with its matrix at every
wavenumber and every input: $\mathcal{Q}(\bm k)|0\rangle$ equals
$\hat{\mat V}(\bm k)$ to $1.9 \times 10^{-14}$, the block of $\mathcal{U}_M$
equals $\tilde{\mat{M}}$ to $3.0 \times 10^{-14}$, the prepared right-hand
side is parallel to $\bm\psi$ to $2.0 \times 10^{-15}$, and the output equals
$q_t(\tilde{\mat{M}})\,\bm\psi/C$, up to the normalization of the right-hand
side, to $5.1 \times 10^{-14}$. Against the finite element solution from a
sparse direct solve, the energy-norm error of the output is
$6.4 \times 10^{-4}$ at $N = 4$ and $7.0 \times 10^{-4}$ at $N = 8$, below
$\epsilon_p = 10^{-3}$, as Equation~\eqref{eq:energy-error} requires with
$\delta = 0$.

The success probability is $2.4 \times 10^{-3}$ at both meshes, the product
of the probability of preparing the right-hand side, $0.070$ at $N = 4$ and
$0.074$ at $N = 8$, and the probability of the inverse, which is limited by
the subnormalization $C = 13.2$ of $q_t$; the latter factor is the expected
cost of QSVT inversion and is removed by $O(C)$ rounds of amplitude
amplification~\cite{lapworth2026probability}. The circuit size, in contrast,
grows with $N$: the CX count of one application of $\mathcal{U}_V$ grows from
$138$ to $310$ when $N$ doubles, most of the increase coming from the two
uniformly controlled rotations, whose $3N^2$ angles make the multiplexor grow
as $N^2$ (Section~\ref{sec:cost}). The simulation is therefore a check that
the gates compose into the inverse, not evidence that this realization
scales; the arithmetic of Appendix~\ref{app:multiplexed} grows only with $F$, and its
logical cost is given next.

\subsection{Resources}
\label{sec:exp-resources}

Table~\ref{tab:resources} gives the logical cost of one application of the
inverse for the square at $N = 1024$ and $\rho = 10$. Each of the $d_t$ queries
to $\mathcal{U}_M$ applies $\mathcal{U}_V$ twice and the controlled material
reflection once, so the total is $d_t$ times the cost of a query in
Section~\ref{sec:cost}.
\begin{itemize}[leftmargin=1.4em,itemsep=2pt,topsep=2pt]
\item \textbf{$\epsilon_p = 10^{-6}$:} $d_t = 145$ queries and $F = 45$, with
$88{,}926$ Toffoli gates and $588$ rotations per query, give
$1.3 \times 10^{7}$ Toffoli gates on $5{,}699$ logical qubits.
\item \textbf{$\epsilon_p = 10^{-3}$:} $d_t = 77$ queries and $F = 34$ give
$4.0 \times 10^{6}$ Toffoli gates on $3{,}576$ logical qubits.
\end{itemize}
Nearly all of the qubits are in the work register of the arithmetic
(Section~\ref{sec:cost}).

\paragraph{Costs not included.}
Three costs are not included in these counts.
\begin{enumerate}[leftmargin=1.4em,itemsep=2pt,topsep=2pt]
\item \textbf{Rotation synthesis:} the single-qubit rotations of
Table~\ref{tab:resources} are counted before synthesis into Clifford and T
gates.
\item \textbf{Phase rotations of QSVT:} the $d_t + 1$ rotations about the input
subspace each mark a subspace of a few qubits, at a number of Toffoli gates that
does not depend on $N$ or $F$.
\item \textbf{Microstructure loading:} the three cells are dyadic, so their
material reflections cost at most $9$ Toffoli gates
(Section~\ref{sec:exp-cells}); a microstructure without a compact
description requires a material oracle whose cost grows with the number of
elements.
\end{enumerate}

\section{Conclusions}
\label{sec:conclusions}

Classical FFT-preconditioned solvers rely on a condition number that does
not grow with the mesh; in a quantum linear solver, this property does not
survive a separately encoded preconditioner. It is recovered by placing the
closed-form inverse square root of the homogeneous operator inside the
factors of the periodic two-phase operator, which gives a block encoding of
the preconditioned operator with subnormalization $1$ and condition number
at most $\rho$ for every $N$. The mesh enters the cost only through the
quantum Fourier transforms and the wavenumber-controlled stage, both
polylogarithmic in $N$, and the condition number, which is the quantity that
sets the number of queries, is a material constant.

Because the condition number is a material constant, the number of
queries depends on the contrast $\rho$ and on the accuracy $\epsilon_p$, but
not on the mesh. For isotropic phases of contrast $\rho$, the degree of the
inverse grows linearly in $\rho$ and does not grow with $N$: at
$\epsilon_p = 10^{-6}$, the degree given by Equation~\eqref{eq:degree} is
$d_t = 145$ at $\rho = 10$ and $d_t \le 1.5 \times 10^{5}$ at $\rho = 10^4$.
Each query applies the fixed-point arithmetic that prepares the vector at
each wavenumber. The precision of this arithmetic grows with $N$, with $d_t$,
and with the subnormalization $C$ (Equation~\eqref{eq:F-choice}); at
$N = 1024$ and $\epsilon_p = 10^{-6}$, one application of the inverse needs
$F = 45$ fractional bits and $1.3 \times 10^{7}$ Toffoli gates. The same
vectors can also be prepared by a table lookup of classically computed
angles, which uses fewer T gates below $N = 4096$ in two dimensions; the cost
of the arithmetic, in contrast, grows with $N$ only through $F$. The
arithmetic needs too many qubits for a state-vector simulation, so in the
simulation of Section~\ref{sec:exp-qiskit}, the rotation at each wavenumber
is built as a multiplexor, whose size grows as $N^2$. With this simulation,
we confirm that the stages compose into the inverse; the arithmetic is
verified separately, one stage at a time.

Three topics are not covered in this paper:
\begin{enumerate}[leftmargin=1.4em,itemsep=2pt,topsep=2pt]
\item \textbf{Mesh and coefficients:} the construction requires a uniform
periodic mesh, a homogeneous reference operator with closed-form
eigenvalues, and isotropic phases. Dirichlet boundaries and general
coefficients are handled by the multilevel construction compared in
Section~\ref{sec:be-M}.
\item \textbf{Gate synthesis and qubit count:} the reported cost excludes
the synthesis of the rotations into Clifford and T gates, and the qubit
count of the reversible arithmetic is not optimized.
\item \textbf{Readout:} the inverse prepares a state proportional to the
solution; reading out the effective conductivity from it is the subject of
future work.
\end{enumerate}
\subsection*{Acknowledgments}
The author acknowledges the Vilas Associate Grant from the University of
Wisconsin Graduate School.

\subsection*{Conflicts of Interest}
The author declares no conflicts of interest.

\subsection*{Data Availability Statement}
The code developed in this work is available as an open-source repository at \url{\repourl}, where the folder \texttt{papers/inversion} contains the scripts that reproduce every table and figure in this paper.

\subsection*{Use of AI Tools}
The author used generative AI assistants, specifically Claude Opus 5.5 (Anthropic), to prepare this manuscript. The author used the tool for code development and text drafting. All content was reviewed, verified, and edited by the author, who takes full responsibility for the accuracy and integrity of the work.

\appendix
\section{Multiplexed Preparation of the Fourier Block}
\label{app:multiplexed}

This appendix constructs the wavenumber-controlled preparation of
Section~\ref{sec:cost}. With the phase $e^{\mathrm{i}(\phi_1+\phi_2)}$ and the
fixed unitary $\mat{T}$ of Equation~\eqref{eq:Vhat-c} separated, what remains
is the unitary
\begin{equation}
\mathcal P = \sum_{\bm k} |\bm k\rangle\langle\bm k| \otimes \mathcal P(\bm k),
\qquad \mathcal P(\bm k)|0\rangle = \bm c(\bm k)/\norm{\bm c(\bm k)}_2 ,
\label{eq:mx-target}
\end{equation}
on the $2n$ qubits of $\bm k$ and two target qubits $t_0$ and $t_1$, for every
$\bm k \ne \bm 0$, with $\bm c = (a, b, \zeta/\sqrt3, \zeta/\sqrt3)$ and $a$,
$b$, and $\zeta$ of Equation~\eqref{eq:abz}; the zero mode is routed as in
Section~\ref{sec:cost}. Its error at $\bm k$ is
$\delta(\bm k) = \norm{\mathcal P(\bm k)|0\rangle - \bm c(\bm k)/\norm{\bm c(\bm k)}_2}_2$.

A single state of this kind, a vector of length $4$, is prepared in standard
state preparation \cite{mottonen2004transformation,shende2006} with three
rotations whose angles are computed when the circuit is built. Multiplexed
over $\bm k$, the angles are tabulated classically for every wavenumber and
applied either by a \emph{multiplexor} of uniformly controlled
rotations~\cite{mottonen2004transformation}, $3N^2$ in all, or by a lookup in
a quantum read-only memory
(QROM)~\cite{babbush2018encoding,low2024trading,berry2019qubitization}
followed by rotations, with $O(N\sqrt b)$ Toffoli gates for angles of $b$
bits. In both, $O(N^2)$ classical operations are needed to tabulate the
angles, and the description of the circuit grows with $N$. The entries $a$,
$b$, and $\zeta$, however, are products of sines and cosines of
$\phi_j = \pi k_j/N$, and they can be evaluated from the bits of $\bm k$ by
arithmetic. This is the structured state preparation of our earlier
work~\citep{suresh2026pyencode}, in which a vector given by a formula is
prepared at a cost polynomial in the number of qubits, applied here to a family
of vectors indexed by $\bm k$. We therefore construct
$\mathcal P = \mathcal A^{\dagger}\mathcal R\,\mathcal A$, where $\mathcal A$
computes $a$, $b$, and $\zeta$ as fixed-point bit strings
(Appendix~\ref{app:mx-part1}), $\mathcal R$ prepares $\bm c/\norm{\bm c}_2$ on
the target qubits from these registers and leaves them unchanged
(Appendix~\ref{app:mx-part2}), and $\mathcal A^{\dagger}$ clears them.
A direct construction of Part 2 would compute $\norm{\bm c}_2$ with a square
root, divide the entries by it, and compute each rotation angle with an inverse
trigonometric function into an \emph{angle register}, a register whose bits
then control the rotations. Each of these functions is a reversible circuit
whose Toffoli count grows as the square of the register width, against linearly
for an addition, and each angle register adds a register of
qubits~\cite{haner2018optimizing}. Sanders et al.~\cite{sanders2019black}
replace the angle register with an inequality test followed by amplitude
amplification; in a multiplexed preparation, however, the number of rounds is
set by the shortest vector, with $\norm{\bm c(\bm k)}_2 \approx \pi/N$ near
$\bm k = \bm 0$, and so grows as $N$. We therefore require the circuit to use
additions, shifts, and controlled rotations only, with no division, square
root, or angle register; its work register must return to $|0\rangle$. The
circuit is constructed in Appendices~\ref{app:mx-part1}
and~\ref{app:mx-part2}, its error is bounded in Appendix~\ref{app:mx-error},
and its gates are counted and compared with the two alternatives in
Appendix~\ref{app:mx-cost}.

\subsection{Part 1: Computing the Bit Strings}
\label{app:mx-part1}

In Part 1, $a$, $b$, and $\zeta$ are computed for every wavenumber at once, as
integers in fixed-point representation and with additions and shifts only. Each
register has $F + 4$ bits, of which $F$ are fractional, so its resolution is
$\Delta = 2^{-F}$. To compute $a$, $b$, and $\zeta$ with shifts and additions,
we first write them, by the sum-to-product identities with
$\theta_\pm = \phi_1 \pm \phi_2$, as
\begin{equation}
\begin{aligned}
a &= \tfrac12(\sin\theta_+ + \sin\theta_-),\\
b &= \tfrac12(\sin\theta_+ - \sin\theta_-),\\
\zeta &= \tfrac12(\cos\theta_- - \cos\theta_+),
\end{aligned}
\label{eq:mx-s2p}
\end{equation}
so the circuit needs the cosine and sine of two angles that are linear in
$\bm k$. These are computed in stages, whose registers and order are shown in
Figure~\ref{fig:mx-part1}.

\begin{figure}[t]
\centering
\resizebox{\linewidth}{!}{%
\begin{quantikz}[row sep={0.8cm,between origins}, column sep=0.35cm, font=\small]
\lstick{$|\bm k\rangle$} & \qwbundle{2n} & & & \gate[3]{\text{add, sub}}\gategroup[7,steps=3,style={dashed,rounded corners,inner xsep=2pt,inner ysep=1pt},label style={label position=above,anchor=north,yshift=0.32cm}]{$\mathcal A$: Part 1} & & & & \gate[7]{\mathcal A^{\dagger}} & \rstick{$|\bm k\rangle$}\\
\lstick{$z_+$: $|0\rangle$} & \qwbundle{F+1} & & & & \gate[5]{\begin{array}{c}\text{rotation mode,}\\ \text{two runs}\end{array}} & & & & \rstick{$|0\rangle$}\\
\lstick{$z_-$: $|0\rangle$} & \qwbundle{F+1} & & & & & & & & \rstick{$|0\rangle$}\\
\lstick{work: $|0\rangle$} & \qwbundle{} & & & & & & & & \rstick{$|0\rangle$}\\
\lstick{$(x_+, y_+)$: $|0\rangle$} & \qwbundle{2(F+4)} & & & & & \gate[3]{\begin{array}{c}\text{sum to}\\ \text{product}\end{array}} & & & \rstick{$|0\rangle$}\\
\lstick{$(x_-, y_-)$: $|0\rangle$} & \qwbundle{2(F+4)} & & & & & & & & \rstick{$|0\rangle$}\\
\lstick{$a, b, \zeta$: $|0\rangle$} & \qwbundle{3(F+4)} & & & & & & \gate[2]{\mathcal R} & & \rstick{$|0\rangle$}\\
\lstick{$t_0, t_1$: $|0\rangle$} & \qwbundle{2} & & & & & & & & \rstick{$\bm c/\norm{\bm c}_2$}
\end{quantikz}}
\caption{The registers and stages of
$\mathcal P = \mathcal A^{\dagger}\mathcal R\,\mathcal A$ for one wavenumber
$\bm k$, read from left to right; the circuit acts on all wavenumbers in
superposition. Part 1, the operation $\mathcal A$ in the dashed box, forms the
angles $z_\pm = (k_1 \pm k_2)/N$ with an adder and a subtractor, runs rotation
mode once for each angle to write
$(x_\pm, y_\pm) \approx (\cos\theta_\pm, \sin\theta_\pm)$, and writes $a$, $b$,
and $\zeta$ by sum to product. The work register holds the range-reduction
flags and the copies that each step keeps for reversibility. Part 2, the
operation $\mathcal R$ of Figure~\ref{fig:mx-circuit}, prepares
$\bm c/\norm{\bm c}_2$ on the two target qubits and leaves $a$, $b$, and
$\zeta$ unchanged; its own work register is not drawn. $\mathcal A^{\dagger}$
then returns every register except $\bm k$ and the target qubits to
$|0\rangle$.}
\label{fig:mx-part1}
\end{figure}

\paragraph{The angles.} The first stage forms the two angles $\theta_\pm$, each
held in units of $\pi$, as $z_\pm = (k_1 \pm k_2)/N$ modulo $2$, on $F + 1$
bits in two's complement, so that the register wraps around at multiples of
$2\pi$. In these units, the value $k_1/N$ needs no Toffoli gates, since
$N = 2^n$: the $n$ bits of $k_1$ are copied with CNOT gates into positions
$F - n$ to $F - 1$. With $k_2/N$ placed in the same way, one adder and one
subtractor then form $z_+$ and $z_-$.

\paragraph{Rotation mode.} From these angles, $\cos\theta$ and $\sin\theta$ are
computed with CORDIC \cite{volder1959cordic} by rotating the vector $(1, 0)$
through the angle $\theta$. The rotation is made in $F$ steps, by the fixed
angles $\mu_i = \arctan 2^{-i}$ for $i = 0, \dots, F - 1$, each in the positive
or the negative direction, so that their signed sum approaches $\theta$. With
$\tan\mu_i = 2^{-i}$, a rotation by one of them needs no multiplication:
\begin{equation}
R(\pm\mu_i) = \cos\mu_i
\begin{bmatrix} 1 & \mp 2^{-i}\\ \pm 2^{-i} & 1 \end{bmatrix},
\label{eq:mx-micro}
\end{equation}
and multiplying a vector $(x, y)$ by the matrix on the right is a shift of the
integers $x$ and $y$ by $i$ bits and two additions. The factor
$\cos\mu_i = (1 + 4^{-i})^{-1/2}$ is the same for both directions, so it is not
applied in each step. The $F$ factors multiply to $1/\Gamma$, where
$\Gamma = \prod_{i<F}(1 + 4^{-i})^{1/2}$ is the \emph{gain}, and the steps
start from $(1/\Gamma, 0)$ instead of $(1, 0)$ to compensate.

\paragraph{The steps.} During the steps, an angle register $z$ holds the part
of $\theta$ that remains to be rotated, starting from $z_+$ or $z_-$, and the
direction of each step is read from the sign bit of $z$: $s_i = +1$ if
$z \ge 0$ and $s_i = -1$ if $z < 0$. Step $i$ maps
\begin{equation}
(x, y) \mapsto (x - s_i 2^{-i} y,\ y + s_i 2^{-i} x),
\qquad z \mapsto z - s_i\tilde\mu_i ,
\label{eq:mx-rotmode}
\end{equation}
where $\tilde\mu_i$ is $\mu_i/\pi$ rounded to a multiple of $\Delta$. Each step
rotates toward $\theta$ by an angle about half that of the previous step, so
the vector may pass $\theta$ but returns toward it in the steps that follow, as
shown in Figure~\ref{fig:mx-rotation} for $\theta = \pi/8$. After the $F$
steps, $(x, y) \approx (\cos\theta, \sin\theta)$, with a remaining angle of at
most $\mu_{F-1} < 2^{1-F}$, of the order of $\Delta$. This remaining angle is
one part of the error; the error as a whole, including the truncation of the
shifts, is bounded in Appendix~\ref{app:mx-error}.

\begin{figure}[H]
\centering
\begin{tikzpicture}[scale=5.5, >=stealth, font=\footnotesize]
\draw[->] (0,0) -- (1.08,0) node[below] {$x$};
\draw[->] (0,0) -- (0,0.72) node[left] {$y$};
\draw[black!40] (1,0) arc[start angle=0, end angle=42, radius=1];
\draw[dashed] (0,0) -- (22.5:1.1) node[right] {$\theta = \pi/8$};
\draw[black!35, ->] (0,0) -- (0.6073,0.0000);
\draw[black!35, ->] (0,0) -- (0.6073,0.6073);
\draw[black!35, ->] (0,0) -- (0.9109,0.3036);
\draw[thick] (0.6073,0.0000) -- (0.6073,0.6073) -- (0.9109,0.3036)
  -- (0.8350,0.5313) -- (0.9014,0.4270) -- (0.9281,0.3706);
\foreach \p/\l/\a in {(0.6073,0.0000)/0/below right, (0.6073,0.6073)/1/above,
  (0.9109,0.3036)/2/below right, (0.8350,0.5313)/3/above left,
  (0.9014,0.4270)/4/above right, (0.9281,0.3706)/5/below right}
  \fill \p circle[radius=0.0065] node[\a] {$\l$};
\node[below left] at (0.6073,0.0000) {$1/\Gamma$};
\end{tikzpicture}
\caption{The first five steps of rotation mode for $\theta = \pi/8$, without
rounding. Point $i$ is the vector $(x, y)$ after $i$ steps, starting from
$(1/\Gamma, 0)$; the directions $s_0, \dots, s_4$ are $+, -, +, -, -$. Each step
rotates by $\mu_i$ toward $\theta$, and the length grows from
$1/\Gamma = 0.607$ toward $1$, the arc of the unit circle. The angle after
step $5$ is $21.8^\circ$, against $22.5^\circ$.}
\label{fig:mx-rotation}
\end{figure}

\paragraph{Range reduction.} The steps reach $\theta$ only if $\theta$ lies
within the largest signed sum of the $\mu_i$, which is $\sum_i\mu_i = 1.74$;
rotation mode therefore converges only for $|\theta| \le 1.74$. The angles
$\theta_\pm$ take every value in $[-\pi, \pi)$, and an angle outside this range
is brought into it by a shift of $\pi$, which only changes the sign of the
result:
\begin{equation}
(\cos\theta, \sin\theta) = -\bigl(\cos(\theta \mp \pi),\ \sin(\theta \mp \pi)\bigr),
\label{eq:mx-range}
\end{equation}
with the upper sign for $\theta \in [\pi/2, \pi)$ and the lower sign for
$\theta \in [-\pi, -\pi/2)$. For an angle in either interval, the circuit
rotates by $\theta \mp \pi$, which lies in $[-\pi/2, \pi/2)$, and starts the
steps from $(-1/\Gamma, 0)$ instead of $(1/\Gamma, 0)$. In the angle
register, this reduction needs no Toffoli gates. The sign bit of $z$ has weight
$-1$, so flipping it changes $z$ by $-1$ if $z \ge 0$ and by $+1$ if $z < 0$,
which is the shift of $\theta$ by $\mp\pi$ modulo $2\pi$. The two leading bits
of $z$ are $01$ for $z \in [1/2, 1)$, $10$ for $z \in [-1, -1/2)$, and $00$ or
$11$ otherwise, so the flag that selects the reduction is their exclusive or,
written with two CNOT gates. Under this flag, the sign bit of $z$ is flipped
with one CNOT gate, and the start value $\pm 1/\Gamma$, a constant, is written
with X and CNOT gates.

\paragraph{Reversibility.} After the reduction, each step overwrites $x$, $y$,
and $z$. To keep the circuit reversible, each step keeps its sign bit and the
shifted copy needed by its second addition until $\mathcal A^{\dagger}$ clears
them \cite{bennett1973logical}; these copies are most of the work register
(Table~\ref{tab:mx-cost}). The two runs, one for $z_+$ and one for $z_-$, give
the four values needed in Equation~\eqref{eq:mx-s2p}.

\paragraph{Sum to product.} From these four values, Equation~\eqref{eq:mx-s2p}
is evaluated with three additions of registers shifted by one bit, since each
factor $\tfrac12$ is a shift of the integer $Z$. Only $a$, $b$, and $\zeta$ are
formed: the factor $1/\sqrt3$ of the entries of $\bm c$ is applied as a
constant in Part 2.

\subsection{Part 2: Computing the Amplitudes}
\label{app:mx-part2}

In Part 2, the bit strings of $a$, $b$, and $\zeta$ are converted into the
amplitudes of the target qubits $t_0$ and $t_1$. For every wavenumber, the
circuit must carry out the map
\begin{equation}
\begin{aligned}
&|a\rangle|b\rangle|\zeta\rangle\,|00\rangle \mapsto
\\
&\quad |a\rangle|b\rangle|\zeta\rangle \otimes
\frac{a|00\rangle + b|01\rangle + \tfrac{\zeta}{\sqrt3}|10\rangle
+ \tfrac{\zeta}{\sqrt3}|11\rangle}{\norm{\bm c}_2},
\end{aligned}
\label{eq:mx-part2}
\end{equation}
where $|a\rangle$ is the basis state of the register that holds the integer of
$a$, and the target qubits are written in the order $|t_0 t_1\rangle$. These
amplitudes are set by rotations of the target qubits, as in the binary-tree
decomposition of standard state preparation: $t_0$ is rotated by an angle $\tau$
that divides the norm between the pair $(a, b)$ and the pair
$(\zeta/\sqrt3, \zeta/\sqrt3)$, and $t_1$ is rotated by an angle $\omega$
within the first pair or by $\pi/4$ within the second. Each of these angles is
the angle of a pair of numbers. With CORDIC vectoring, the angle of a pair is
obtained as a sequence of decision bits, one per step, and these bits control
fixed rotations of a target qubit directly. The arithmetic and the rotation
stage are shown in Figure~\ref{fig:mx-circuit}.

\begin{figure*}[t]
\centering
\begin{minipage}[b]{0.38\linewidth}
\centering
\begin{tikzpicture}[
  font=\small, >=stealth,
  reg/.style={draw, minimum width=0.9cm, minimum height=0.5cm},
  vec/.style={draw, rounded corners=2pt, minimum width=1.7cm, minimum height=0.55cm, fill=black!5},
  bits/.style={font=\footnotesize, inner sep=1pt}]
\node[reg] (v0) at (0,0) {$a$};
\node[reg] (v1) at (1.1,0) {$b$};
\node[reg] (v2) at (2.95,0) {$\zeta$};
\node[vec] (n10) at (0.55,1.3) {vectoring};
\node[vec] (n11) at (2.95,1.3) {$\times\sqrt{2/3}\,\Gamma$};
\node[vec] (n00) at (1.75,2.75) {vectoring};
\draw[->] (v0) -- (n10); \draw[->] (v1) -- (n10);
\draw[->] (v2) -- (n11);
\draw[->] (n10) -- node[left, bits] {$\Gamma r$} (n00);
\draw[->] (n11) -- node[right, bits] {$\sqrt{2/3}\,\Gamma\zeta$} (n00);
\draw[->, dashed] (n00) -- ++(0,0.75) node[above, bits] {$\Gamma^2\norm{\bm c}_2$ (uncomputed)};
\node[bits, left=2pt of n10, align=right] {sign of $a$,\\ $\omega$ bits};
\node[bits, right=2pt of n00] {$\tau$ bits};
\node[font=\footnotesize] at (1.5,-0.65) {(a)};
\end{tikzpicture}
\end{minipage}\hfill
\begin{minipage}[b]{0.60\linewidth}
\centering
\scalebox{0.78}{\begin{quantikz}[row sep={0.55cm,between origins}, column sep=0.1cm, font=\small]
\lstick{$t_0$} & \gate{R(\mu_\Sigma)} & \gate{R(-2\mu_i)} & \octrl{4} & & & \ctrl{1} & & \ctrl{3} & \\
\lstick{$t_1$} & & & & \gate{R(\mu_\Sigma)} & \gate{R(-2\mu_i)} & \gate{R(\tfrac{\pi}{4} - \mu_\Sigma)} & & & \\
\lstick{$\tau$ bits} & \qwbundle{F} & \ctrl{-2} & & & & & & & \\
\lstick{$\omega$ bits, sign of $a$} & \qwbundle{F+1} & & \ctrl{1} & & & & & \gate[2]{\text{uncompute}} & \\
\lstick{temporary: $|0\rangle$} & \qwbundle{F+1} & & \targ{} & & \ctrl{-3} & & \gate{Z} & & \rstick{$|0\rangle$}
\end{quantikz}}\\[-0.2em]
{\footnotesize (b)}
\end{minipage}
\caption{The circuit $\mathcal R$ for $\bm c = (a, b, \zeta/\sqrt3,
\zeta/\sqrt3)$. Panel (a) shows the arithmetic; read the arrows upward. The
value $\Gamma r$, with $r = (a^2 + b^2)^{1/2}$, is passed from the vectoring of
$(a, b)$ to the root; no vectoring is needed for the two equal hourglass
entries, and $\sqrt{2/3}\,\Gamma\zeta$, with the same gain, is formed by a
constant product. Only the decision bits of $\omega$ and $\tau$ and the sign
bit of $a$ are passed to the rotation stage. Panel (b) shows the rotation
stage. Each gate with index $i$ represents $F$ gates, one per decision bit.
Qubit $t_0$ is rotated by the $\tau$ bits. The Toffoli column copies the
$\omega$ bits and the sign bit of $a$ into temporary qubits when $t_0 = 0$, so
that the rotations of $t_1$, and the $Z$ gate that restores the sign of $a$,
act only on that branch; the rotation $R(\pi/4 - \mu_\Sigma)$ controlled by
$t_0$ completes the fixed angle $\pi/4$ of the hourglass node. The temporary
qubits are uncomputed by measurement with a phase correction controlled by
$t_0$.}
\label{fig:mx-circuit}
\end{figure*}

\paragraph{One pair.} Consider a pair $(x_0, y_0)$ of register values, not
normalized and of either sign, with norm $L = (x_0^2 + y_0^2)^{1/2} > 0$. The
state to prepare on one qubit is the normalized pair,
\begin{equation}
\frac{x_0}{L}\,|0\rangle + \frac{y_0}{L}\,|1\rangle = R(\theta)|0\rangle,
\qquad \theta = \operatorname{atan2}(y_0, x_0),
\label{eq:mx-pair}
\end{equation}
where $R(\theta) = R_y(2\theta)$, whose matrix is the plane rotation by
$\theta$ of Equation~\eqref{eq:mx-micro}. The state depends only on the angle
$\theta$ of the pair, not on its norm $L$, so the circuit needs to compute an
angle and apply it as a rotation, and it does not need to divide by $L$.

\paragraph{Vectoring.} This angle is computed by CORDIC in vectoring mode
\cite{volder1959cordic}, with the steps of Equation~\eqref{eq:mx-micro}; the
direction of each step is now chosen from the sign of $y$, so that the steps
rotate the pair onto the positive $x$ axis. The \emph{decision bit} $d_i$ of
step $i$ is the sign bit of $y$, equal to $1$ if $y < 0$ and $0$ otherwise.
With the direction $s_i = 1 - 2d_i$, step $i$ maps
\begin{equation}
(x, y) \mapsto (x + s_i 2^{-i} y,\ y - s_i 2^{-i} x),
\label{eq:mx-vecmode}
\end{equation}
with arithmetic (floor) shifts of the integers $x$ and $y$. Each step rotates
the pair by $-s_i\mu_i$, so after $F$ steps the pair has been rotated by
$-\theta_d$, where
\begin{equation}
\theta_d = \sum_{i<F} (1 - 2d_i)\,\mu_i = \mu_\Sigma - 2\sum_{i<F} d_i\,\mu_i ,
\qquad \mu_\Sigma = \sum_{i<F}\mu_i .
\label{eq:mx-theta}
\end{equation}
The rotated pair lies near the positive $x$ axis, so $\theta_d$ approximates
$\theta$. This holds only within the convergence range $\sum_i\mu_i = 1.74$,
which covers $|\theta| \le \pi/2$ but not the left half-plane, so $x_0 \ge 0$
is required; this is arranged in the tree below. By
Equation~\eqref{eq:mx-theta}, the decision bits encode $\theta_d$, one bit per
step, so no register holds the angle; and since each step multiplies the
length of the pair by $(1 + 4^{-i})^{1/2}$, the final $x$ register holds
$x_F \approx \Gamma L$, the norm times the gain, without a square root.

\paragraph{From decision bits to a rotation.} The angle $\theta_d$ is applied
as a rotation without being computed, since rotations about one axis add,
$R(\alpha)R(\beta) = R(\alpha + \beta)$; by Equation~\eqref{eq:mx-theta},
\begin{equation}
R(\theta_d) = R(\mu_\Sigma)\prod_{i<F} R(-2\mu_i)^{d_i} :
\label{eq:mx-ctrlrot}
\end{equation}
the fixed rotation $R(\mu_\Sigma)$ followed by the rotations $R(-2\mu_i)$, each
controlled by its decision bit, as in Burge et al.~\cite{burge2024quantum}. The
angles $\mu_\Sigma$ and $\mu_i$ are known when the circuit is built, so the
circuit applies $R(\theta_d)$ with no arithmetic on $\theta_d$. The decision
bits must therefore remain available until the rotations are applied: each step
keeps its bit $d_i$ and the shifted copy of $x$ until the uncomputation
\cite{bennett1973logical}.

\paragraph{The tree.} A vector of length $4$ is normalized by normalizing pairs
at two levels, as in the binary-tree decomposition of standard state
preparation \cite{grover2002creating,mottonen2004transformation}. The pair
$(a, b)$ has norm $r = (a^2 + b^2)^{1/2}$, and the hourglass pair
$(\zeta/\sqrt3, \zeta/\sqrt3)$ has norm $\sqrt{2/3}\,\zeta$. In terms of these
norms, with the first factor of each term on $t_0$ and the second on $t_1$,
\begin{equation}
\begin{gathered}
\begin{aligned}
\frac{\bm c}{\norm{\bm c}_2} = {}&\cos\tau\,|0\rangle
\bigl(\cos\omega\,|0\rangle + \sin\omega\,|1\rangle\bigr)
\\
&+ \sin\tau\,|1\rangle\,\frac{|0\rangle + |1\rangle}{\sqrt2},
\end{aligned}
\\
\omega = \operatorname{atan2}(b, a), \quad
\tau = \operatorname{atan2}\bigl(\sqrt{2/3}\,\zeta,\ r\bigr).
\end{gathered}
\label{eq:mx-tree}
\end{equation}
The hourglass pair has the angle $\pi/4$ for every wavenumber, so only $\omega$
and $\tau$ are computed, node by node (Figure~\ref{fig:mx-circuit}(a)).
\begin{itemize}[leftmargin=1.4em,itemsep=2pt,topsep=2pt]
\item \textbf{The pair $(a, b)$:} a vectoring gives $\omega$ in its decision
bits and $\Gamma r$ in its $x$ register. Since vectoring requires $x_0 \ge 0$,
the pair is first negated by bitwise complement when $a < 0$, with CNOT gates
controlled by the sign bit of $a$; the negation changes $\omega$ by $\pi$, and
the rotation stage restores the sign.
\item \textbf{The hourglass pair:} its angle is known, so no vectoring is
needed. Its norm must enter the root with the same gain as $\Gamma r$, so the
circuit forms the constant product $\sqrt{2/3}\,\Gamma\zeta$ by shifts and
additions over the $n_s$ signed digits of $\sqrt{2/3}\,\Gamma$, the leading
digit a copy. Since $\zeta \ge 0$, no negation is needed.
\item \textbf{The root:} a second vectoring, of
$(\Gamma r, \sqrt{2/3}\,\Gamma\zeta)$, gives $\tau$ in its decision bits, since
a common factor does not change the angle of a pair. Its first entry is
nonnegative, so no negation is needed.
\end{itemize}
The two vectorings run in place on the registers of $a$ and $b$, and the
constant product writes $\sqrt{2/3}\,\Gamma\zeta$ into a new register of
$F + 4$ bits; all three computations are reversed after the rotation stage. The
norms $\Gamma r$ and $\Gamma^2\norm{\bm c}_2$ are never read.

\paragraph{The rotation stage.} At each level of the tree, one target qubit is
rotated, as in the uniformly controlled rotations of M\"ott\"onen et
al.~\cite{mottonen2004transformation}, with the angles taken from decision bits
instead of a table (Figure~\ref{fig:mx-circuit}(b)). Qubit $t_0$ is rotated by
$R(\tau)$, applied as in Equation~\eqref{eq:mx-ctrlrot} with the decision bits
of $\tau$ as controls. Qubit $t_1$ must then be rotated by $\omega$ when
$t_0 = 0$ and by $\pi/4$ when $t_0 = 1$. For this, $t_1$ receives the fixed
rotation $R(\mu_\Sigma)$, the rotations $R(-2\mu_i)$, each controlled both by
its decision bit $d_i$ of $\omega$ and by $t_0 = 0$, and the rotation
$R(\pi/4 - \mu_\Sigma)$ controlled by $t_0$, which give
\begin{equation}
\begin{aligned}
t_0 = 0:&\quad R(\mu_\Sigma)\prod_{i<F} R(-2\mu_i)^{d_i} = R(\theta_d) \approx R(\omega),
\\
t_0 = 1:&\quad R(\mu_\Sigma)\,R(\tfrac{\pi}{4} - \mu_\Sigma) = R(\tfrac{\pi}{4}).
\end{aligned}
\label{eq:mx-stage}
\end{equation}
The sign of $a$, removed by the negation, is restored on the same branch: the
branch $t_0 = 0$ is multiplied by $-1$ when $a < 0$, since
$R(\theta + \pi) = -R(\theta)$. Each of these operations has two conditions,
and one Toffoli gate writes their logical AND into a temporary qubit, which then
controls the rotation or, for the sign, receives a $Z$ gate. Each temporary
qubit is uncomputed by measurement, with a phase correction by a CZ gate on its
two conditions \cite{babbush2018encoding}. The stage uses $2F + 1$ controlled
rotations, $2F$ of them controlled by decision bits and one,
$R(\pi/4 - \mu_\Sigma)$, controlled by $t_0$, two uncontrolled rotations,
$F + 3$ Toffoli gates, and $F + 2$ temporary qubits.

\paragraph{Example.} For the wavenumber of Example~\ref{ex:precision},
$\bm k = (3, 1)$ at $N = 16$ and $F = 12$, so that $\theta_+ = \pi/4$,
$\theta_- = \pi/8$, and $\Delta = 2^{-12}$, rotation mode writes
$(\cos\theta_+, \sin\theta_+) = (0.70752, 0.70605)$ and
$(\cos\theta_-, \sin\theta_-) = (0.92334, 0.38257)$, in error by $4.6\Delta$
and $2.3\Delta$, and sum to product gives the entries $a$, $b$, and $\zeta$ of
Example~\ref{ex:precision}. Since $a > 0$, the pair $(a, b)$ is vectored
without negation; its vectoring ends with $\Gamma r = 0.93506$, against the
exact $0.93623$, and the constant product is
$\sqrt{2/3}\,\Gamma\zeta = 0.14526$, against $0.14573$. The two vectorings
leave the decision bits $d_0 \dots d_{11}$
\begin{equation}
\begin{aligned}
\omega:&\quad \texttt{011000100101}, && \text{encoded angle } 0.28942,\\
\tau:&\quad \texttt{011010110000}, && \text{encoded angle } 0.15384,
\end{aligned}
\label{eq:mx-example2}
\end{equation}
where each encoded angle is $\mu_\Sigma - 2\sum_i d_i\mu_i$ by
Equation~\eqref{eq:mx-theta}, against the exact $\omega = 0.28934$ and
$\tau = 0.15442$: errors of $0.3\Delta$ and $2.4\Delta$. With these bits as
controls, the rotation stage prepares the vector of Example~\ref{ex:precision},
with $\delta(\bm k) = 5.9 \times 10^{-4} = 2.4\Delta$.

\subsection{Error}
\label{app:mx-error}

The error $\delta(\bm k)$ arises in three places. First, rotation mode in
Part 1 is stopped after $F$ steps, with rounded angle constants and truncated
shifts, so $a$, $b$, and $\zeta$ carry an absolute error of a few units of
$\Delta$. Second, the entries are normalized in Part 2, so this absolute error
becomes an error relative to $\norm{\bm c}_2$, and it is magnified for a short
vector: at $\bm k = (1, 0)$ and $N = 1024$,
$\norm{\bm c}_2 = \sin(\pi/1024) = 3.1 \times 10^{-3}$, and an absolute error
of $\Delta$ in $a$ is a relative error of about $330\Delta$. Third, the shifts
of each vectoring in Part 2 are truncated, and the resulting error is passed up
the tree with the norms that the vectoring outputs.

\begin{proposition}[Error of the multiplexed preparation]
\label{prop:mx-error}
Let $N = 2^n \ge 4$ and $F \ge n$. There is a constant $C_\delta$, independent
of $N$, $F$, and $\bm k$, such that for every $\bm k \ne \bm 0$
\begin{equation}
\delta(\bm k) \le \frac{C_\delta\,F\,2^{-F}}{\norm{\bm c(\bm k)}_2},
\qquad \norm{\bm c(\bm k)}_2 \ge \sin(\pi/N) \ge \frac{2}{N},
\label{eq:mx-error}
\end{equation}
so that $\delta(\bm k) \le \tfrac12\,C_\delta\,F N\,2^{-F}$.
\end{proposition}

\begin{proof}
The argument has five steps.
\begin{itemize}[leftmargin=1.4em,itemsep=2pt,topsep=2pt]
\item \textbf{Rotation mode:} after $F$ steps the remaining angle is at most
$\mu_{F-1} < 2^{1-F}$, the error from rounding the $F$ angle constants is at
most $F/2$ units of the angle register, and the truncation of the shifts in
each step is at most $\sqrt2\,\Delta$ and is enlarged by at most
$\prod_{j \ge 2}(1 + 4^{-j})^{1/2} = 1.0415$ in the steps that follow. The
computed $(\cos\theta_\pm, \sin\theta_\pm)$ therefore have an error of order
$F\Delta$ \cite{hu1992quantization}.
\item \textbf{Sum to product:} by Equation~\eqref{eq:mx-s2p}, the entries are
half sums and half differences of these values, with a truncation of at most
$\Delta/2$ per shift, so the computed vector $\tilde{\bm c}$ satisfies
$\norm{\tilde{\bm c} - \bm c}_2 \le \delta_{\mathcal A}$, with
$\delta_{\mathcal A}$ of order $F\Delta$.
\item \textbf{One vectoring:} for a pair of norm $L \le 2$, which covers the
inputs of the root, scaled by $\Gamma$, the error of the direction given by the
decision bits, weighted by $L$, and the error of the length $x_F/\Gamma$ are
each bounded by the same argument by a constant times $2^{-F}L + F\Delta$
\cite{hu1992quantization}, so by a constant times $F\Delta$.
\item \textbf{The tree:} for a node $\nu$ whose pair has exact direction
$\bm u_\nu$ and norm $L_\nu$, and prepared direction $\tilde{\bm u}_\nu$, let
$\varepsilon_\nu = L_\nu\norm{\tilde{\bm u}_\nu - \bm u_\nu}_2$. Replacing the
exact angles by the computed ones one node at a time, from the leaves to the
root, changes the prepared vector by $\varepsilon_\nu/\norm{\tilde{\bm c}}_2$
at node $\nu$, since its ancestors are exact and the rotations below it act
unitarily, so the prepared vector differs from
$\tilde{\bm c}/\norm{\tilde{\bm c}}_2$ by at most
$\sum_\nu \varepsilon_\nu/\norm{\tilde{\bm c}}_2$. The hourglass node has an
exact angle. The inputs of the root are the lengths of its two children, each
with an error of order $F\Delta$, and replacing a pair $\bm p$ by $\bm p'$
changes its direction, weighted by its norm, by at most
$2\norm{\bm p' - \bm p}_2$. Every $\varepsilon_\nu$ is therefore of order
$F\Delta$.
\item \textbf{Combination:} by the same inequality,
$\norm{\tilde{\bm c}/\norm{\tilde{\bm c}}_2 - \bm c/\norm{\bm c}_2}_2
\le 2\delta_{\mathcal A}/\norm{\bm c}_2$. If
$\delta_{\mathcal A} \le \norm{\bm c}_2/2$, then
$\norm{\tilde{\bm c}}_2 \ge \norm{\bm c}_2/2$, and the two parts add to a
constant times $F\Delta/\norm{\bm c}_2$. Otherwise $\delta(\bm k) \le 2 <
4\delta_{\mathcal A}/\norm{\bm c}_2$ holds because both vectors are unit
vectors. For the lower bound, let $\gamma_j = \sin^2\phi_j$, so that
$\norm{\bm c}_2^2 = \gamma_1 + \gamma_2 - \tfrac43 \gamma_1 \gamma_2$, which is
linear in each $\gamma_j \in [0, 1]$. For $\bm k \ne \bm 0$, one of them, say
$\gamma_1$, is at least $\sin^2(\pi/N)$, and at $\gamma_2 = 0$ and
$\gamma_2 = 1$ the value is $\gamma_1$ and $1 - \gamma_1/3 \ge 2/3$. Hence
$\norm{\bm c}_2 \ge \sin(\pi/N)$, and $\sin(\pi/N) \ge 2/N$.\qedhere
\end{itemize}
\end{proof}

The factor $F$ in Equation~\eqref{eq:mx-error} arises because the truncation
errors of the $F$ steps of each CORDIC run add \cite{hu1992quantization}, and
the factor $N$ because the shortest vectors have norm $\sin(\pi/N)$.
Proposition~\ref{prop:mx-error} gives the form of the error, not a value of
$C_\delta$ from which $F$ can be chosen.

\begin{figure}[H]
\centering
\includegraphics[width=\linewidth]{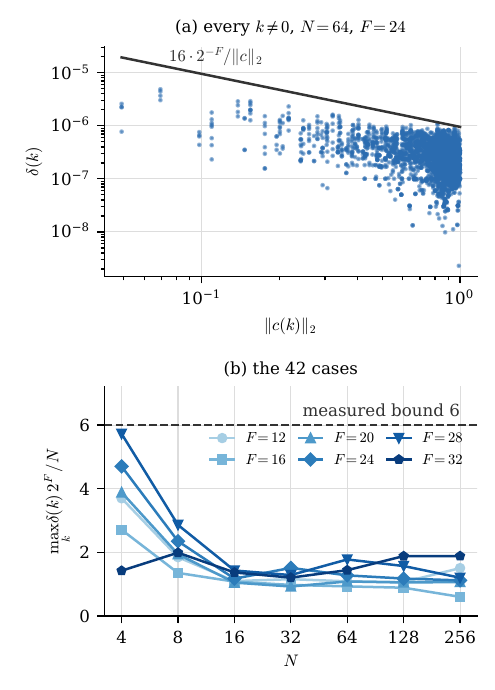}
\caption{The measured error $\delta(\bm k)$ of the multiplexed preparation,
simulated bit by bit. (a) Every $\bm k \ne \bm 0$ at $N = 64$ and $F = 24$,
against $\norm{\bm c(\bm k)}_2$; every point lies below the line
$16 \cdot 2^{-F}/\norm{\bm c}_2$, and the points at the left are the short
vectors near $\bm k = \bm 0$. (b) The largest error over $\bm k$, divided by
$N\,2^{-F}$, for $N = 4$ to $256$ and $F = 12$ to $32$. Every value lies below
the measured bound $6$; the values near $6$ occur only at $N = 4$ and $8$. No
bound is implied for configurations outside these $42$ cases.}
\label{fig:mx-error}
\end{figure}

\paragraph{Measured error.} To obtain a bound with a numerical constant, we
simulate Parts 1 and 2 bit by bit for every $\bm k \ne \bm 0$, for $N = 4$ to
$256$ and $F = 12, 16, \dots, 32$, $42$ cases in all; in every case the work
register returns to zero. The error at every wavenumber for $N = 64$ and
$F = 24$ is shown in Figure~\ref{fig:mx-error}(a). It follows the dependence on
$1/\norm{\bm c}_2$ of Proposition~\ref{prop:mx-error}, with
$\delta(\bm k)\,\norm{\bm c(\bm k)}_2 \le 16 \cdot 2^{-F}$ at every $\bm k$.
Over all $42$ cases, the largest error at each $N$ and $F$, divided by
$N\,2^{-F}$, is shown in Figure~\ref{fig:mx-error}(b). This ratio is at most
$5.73$, reached at $N = 4$, and at most $1.89$ for $N \ge 16$, and it does not
grow systematically with $F$ over this range, although a factor $F$ is allowed
by Proposition~\ref{prop:mx-error}. From these cases we take the measured bound
$6N\,2^{-F}$ of Equation~\eqref{eq:delta-bound}, an observation, not a theorem;
for $N \ge 16$ it exceeds the measured error by a factor of about $3$.

\subsection{Cost and Comparison}
\label{app:mx-cost}

The gates of the circuit are mostly Toffoli gates in additions, which we count
with the adders of Gidney \cite{gidney2018halving}, at $F + 3$ Toffoli gates
per addition of two $(F+4)$-bit integers, for every stage computed and
uncomputed. Each step of rotation mode consists of two additions and one update
of the angle register, $2(F + 3) + F$ Toffoli gates, and each step of vectoring
of two additions, $2(F + 3)$. Table~\ref{tab:mx-cost} lists the counts of the
stages at $N = 1024$ and $F = 45$, the setting of Section~\ref{sec:cost}. The
two runs of rotation mode cost more than the two vectorings, and no count
depends on $N$; in the repository, the counts are checked against the gate
lists of the circuit.

\begin{table}[H]
\caption{Toffoli gates of the multiplexed preparation $\mathcal P$ at
$N = 1024$ and $F = 45$, by stage, compute plus uncompute. Most of the gates
are in rotation mode. The work register holds $5{,}672$ qubits: $5{,}625$ for
the arithmetic, mostly the shifted copies that each step keeps until the
uncomputation, and the $47$ temporary qubits of the rotation stage. The $20$
qubits of the wavenumber and the two target qubits are not included, nor are
the $38$ Toffoli gates that route $\bm k = \bm 0$ (Section~\ref{sec:cost}).}
\label{tab:mx-cost}
\centering\small
\begin{tabular}{llr}
\toprule
part & stage & Toffoli\\
\midrule
1 & angles $\theta_\pm$: adder and subtractor & $180$\\
1 & cosines and sines: two rotation-mode runs & $25{,}380$\\
1 & $a$, $b$, $\zeta$: sum to product & $288$\\
2 & $\omega$: vectoring of $(a, b)$ & $8{,}640$\\
2 & $\sqrt{2/3}\,\Gamma\zeta$: constant product, $n_s = 14$ & $1{,}248$\\
2 & $\tau$: vectoring of $(\Gamma r, \sqrt{2/3}\,\Gamma\zeta)$ & $8{,}640$\\
2 & rotation stage: selection of the $\omega$ bits & $48$\\
\midrule
 & total & $44{,}424$\\
\bottomrule
\end{tabular}

\end{table}

Summing the stages for general $F$ gives Proposition~\ref{prop:mx-main}, the
result used in Section~\ref{sec:cost}.

\begin{proposition}[Multiplexed preparation]
\label{prop:mx-main}
Let $N = 2^n \ge 4$ and $F \ge n$. The circuit of
Appendices~\ref{app:mx-part1} and~\ref{app:mx-part2}, with $F$ fractional bits,
implements Equation~\eqref{eq:mx-target} with a work register of
$2F^2 + 35F + 47$ qubits that returns to $|0\rangle$. It uses
\begin{equation}
20F^2 + 59F + 21 + 2(F + 3)(n_s - 1)
\label{eq:mx-toffoli}
\end{equation}
Toffoli gates, where $n_s$ is the number of signed digits of
$\sqrt{2/3}\,\Gamma$ on $F$ bits ($n_s = 14$ at $F = 45$), $2F + 1$ controlled
single-qubit rotations, $2F$ of them controlled by decision bits and one by
$t_0$, and $2$ uncontrolled ones; no count depends on $N$. Its error satisfies
$\delta(\bm k) \le \tfrac12\,C_\delta\,F N\,2^{-F}$ for every
$\bm k \ne \bm 0$ (Proposition~\ref{prop:mx-error}).
\end{proposition}

\paragraph{Comparison with tabulated angles.} We compare the circuit with the
multiplexor and the table lookup at a common accuracy $\delta$, a bound on
$\max_{\bm k}\delta(\bm k)$, in Table~\ref{tab:mx-options}. For the circuit,
we take $F = \lceil \log_2(6N/\delta) \rceil$ from the measured bound of
Equation~\eqref{eq:delta-bound}, an observation for $N \le 256$ that we
extrapolate to larger $N$. For the lookup,
$b = \lceil \log_2(2\pi/\delta) \rceil$ bits per angle suffice, about
$\log_2 N$ fewer than $F$, since the tabulated angles are exact before rounding
and are not affected by the magnification of short vectors. The lookup is
costed with the formulas of Berry et al.~\cite[Sec.~3.2]{berry2019qubitization}
for clean ancillas, with about $3b$ further Toffoli gates to apply its rotations
by phase-gradient addition~\cite{gidney2018halving}. T gates are counted as $4$
per Toffoli gate~\cite{jones2013low} and $3\log_2(1/\epsilon)$ per rotation
synthesized to accuracy $\epsilon$~\cite{ross2016optimal}, with the synthesis
errors together at most $\delta$; a controlled rotation is written as two
single-qubit rotations and two CNOT gates, and so counts as two rotations. We
take $\delta = 1.96 \times 10^{-10}$, for which $b = 35$ and $F = 45$ at
$N = 1024$.

\begin{table*}[t]
\caption{Estimated cost of one multiplexed preparation of
$\bm c(\bm k)/\norm{\bm c(\bm k)}_2$ for the three constructions, at the
accuracy $\delta = 1.96 \times 10^{-10}$, with $F$ from the measured bound of
Equation~\eqref{eq:delta-bound}, extrapolated beyond $N = 256$. Read each row
across: the construction with the smallest T count is the cheapest at that $N$.
The T counts include the synthesis of the rotations, which the Toffoli counts
and Table~\ref{tab:mx-cost} do not. The circuit of this appendix is counted
gate by gate; the multiplexor and the table lookup are estimated from published
cost formulas, not built. No Toffoli gates or work qubits are used in the
multiplexor, so it is not listed in the last four columns. The table is written
by \texttt{makeTables.py options} in the folder \texttt{papers/inversion} of
the repository.}
\label{tab:mx-options}
\centering\footnotesize
\setlength{\tabcolsep}{4pt}
\begin{tabular}{rrrrrrrrr}
\toprule
 & & \multicolumn{3}{c}{T gates} & \multicolumn{2}{c}{Toffoli gates} & \multicolumn{2}{c}{work qubits}\\
\cmidrule(lr){3-5}\cmidrule(lr){6-7}\cmidrule(lr){8-9}
$N$ & $F$ & multiplexor & table lookup & this circuit & table lookup & this circuit & table lookup & this circuit\\
\midrule
$16$ & $39$ & $9.6 \times 10^{4}$ & $1.3 \times 10^{3}$ & $1.5 \times 10^{5}$ & $335$ & $33{,}666$ & $210$ & $4{,}454$\\
$64$ & $41$ & $1.7 \times 10^{6}$ & $4.9 \times 10^{3}$ & $1.7 \times 10^{5}$ & $1{,}235$ & $37{,}116$ & $630$ & $4{,}844$\\
$256$ & $43$ & $2.9 \times 10^{7}$ & $1.9 \times 10^{4}$ & $1.8 \times 10^{5}$ & $4{,}835$ & $40{,}734$ & $2{,}310$ & $5{,}250$\\
$1024$ & $45$ & $5.1 \times 10^{8}$ & $7.7 \times 10^{4}$ & $2.0 \times 10^{5}$ & $19{,}235$ & $44{,}424$ & $9{,}030$ & $5{,}672$\\
$2048$ & $46$ & $2.1 \times 10^{9}$ & $1.5 \times 10^{5}$ & $2.1 \times 10^{5}$ & $38{,}435$ & $46{,}427$ & $17{,}990$ & $5{,}889$\\
$4096$ & $47$ & $8.7 \times 10^{9}$ & $3.1 \times 10^{5}$ & $2.2 \times 10^{5}$ & $76{,}835$ & $48{,}374$ & $35{,}910$ & $6{,}110$\\
\bottomrule
\end{tabular}

\end{table*}

In Table~\ref{tab:mx-options}, the multiplexor, whose T count grows as $N^2$
with its $3N^2$ rotations, has fewer T gates than the circuit only up to
$N = 16$. The table lookup has the fewest T gates at every $N$ up to $2048$; at
$N = 1024$, its Toffoli count is $2.3$ times smaller than that of the circuit,
but it uses $9{,}030$ work qubits, against $5{,}672$, and its Toffoli count
grows as $N$, or as $N^{3/2}$ on a three-dimensional mesh. The circuit, whose
Toffoli count grows with $N$ only through $F$, has the fewest T gates from
$N = 4096$. This crossover rests on two estimates: the published cost formulas
of the lookup, and the measured bound $6N\,2^{-F}$ extrapolated beyond
$N = 256$. Since the measured error for $N \ge 16$ is at most
$1.89\,N\,2^{-F}$, the extrapolation is more likely to overstate the cost of
the circuit than to understate it. The extra $\log_2 N$ bits of the circuit are
due to the few short vectors near $\bm k = \bm 0$, so a combination of the two,
with these vectors prepared from a table and the rest by the circuit, would need
fewer bits; it has not been constructed.

\bibliographystyle{unsrt}
\bibliography{poissonInversionRef}

\end{document}